\documentstyle[11pt]{article}
\newtheorem {th}{Theorem}[section]
\newtheorem {conj}[th]{Conjecture}
\newtheorem {pblm}[th]{Problem}
\newtheorem {prop}[th]{Proposition}
\newtheorem {pr}[th]{Proposition}
\newtheorem {lem}[th]{Lemma}
\newtheorem {cor}[th]{Corollary}
\newtheorem {eg}{Example}
\newtheorem {defn}{Definition}

\newcommand{\lab}{\label}

\newcommand{\ben}{\begin{enumerate}}
\newcommand{\een}{\end{enumerate}}

\newcommand{\bea}{\begin{eqnarray}}
\newcommand{\ba}{\begin{array}}
\newcommand{\bean}{\begin{eqnarray*}}

\newcommand{\ea}{\end{array}}
\newcommand{\eea}{\end{eqnarray}}
\newcommand{\eean}{\end{eqnarray*}}
\newcommand{\be}{\begin{equation}}
\newcommand{\ee}{\end{equation}}

\def\qed{\hfill \Box}

\def\f{\psi}

\def\eps{\epsilon}

\newcommand{\Z}{Z\!\!\!Z}
\newcommand{\bbbz}{Z\!\!\!Z}

\def\P{{\bf P}}
\def\p{\delta}
\def\dd{\delta}

\def\R{\hbox{I\kern-.2em\hbox{R}}}

\def\pt{phase transition }

\def\calrdJ{\rho^d(J)}
\def\calk{{\cal K}}

\def\one{{\bf 1}}

\def\fb{\overline{f}}
\def\hb{\overline{h}}
\def\K{{\cal K}}

\def\op{{\bf Op}_J} 
\def\PP{{\cal P}}
\def\ssp{\langle\PP_+\rangle}
\def\sspj{\langle\PP_+(J)\rangle}
\def\lts{L^2 (\S)}
\def\ltsg{L^2 (\S/\zero)}

\def\zero{{\hat{0}}}

\def\FF(#1,#2,#3){{\{#1 \stackrel{#3}{\leftrightarrow} #2 \}}}
\def\F(#1,#2){{\{\rho \stackrel{#2}{\leftrightarrow} #1 \}}}
\def\FFFF(#1,#2,#3,#4){{\{\rho \stackrel{#2}{\leftrightarrow} #1,
\rho \stackrel{#4}{\leftrightarrow} #3 \}}}
\def\FFFFC(#1,#2,#3,#4){{\rho \stackrel{#2}{\leftrightarrow} #1,
\rho \stackrel{#4}{\leftrightarrow} #3 }}
\def\Fi(#1,#2){{\{\rho_i \stackrel{#2}{\leftrightarrow} #1 \}}}
\def\FFC(#1,#2,#3){{#1 \stackrel{#3}{\leftrightarrow} #2 }}
\def\FC(#1,#2){{ \rho \stackrel{#2}{\leftrightarrow} #1 }}
\def\FCi(#1,#2){{ \rho_i \stackrel{#2}{\leftrightarrow} #1 }}
\def\FCpi(#1,#2){( \rho_i \stackrel{#2}{\leftrightarrow} #1 )}
\def\FS(#1){{ \rho {\leftrightarrow} #1 }}
\def\FCbar(#1,#2){{ \rho {\leftrightarrow} #1 \mbox{ \rm in } \bar\eta_{#2}}}

\def\mufree{\mu^{\rm free}}
\def\J{{\cal J}}

\def\KJ{{\cal K}_J}
\def\KJP{{\cal K}_{J'}}

\def\S{{\bf S}}
\def\dx{d{\bf x}}

\def\poi{{\bigodot}}

\begin{document}

\begin{titlepage}

\title{Robust Phase Transitions for Heisenberg and Other Models 
on General Trees}
\author{Robin Pemantle\thanks{Research partially
supported by a Presidential Faculty Fellowship and a Sloan Foundation
Fellowship.}
\and Jeffrey E. Steif\thanks{Research supported by
grants from the Swedish Natural Science Research Council and from
the Royal Swedish Academy of Sciences.}
 \and \\
\it  University of Wisconsin-Madison
   and Chalmers University of Technology } 
\date{}
\maketitle

\begin{abstract}
We study several statistical mechanical models on a general tree.
Particular attention is devoted to the classical
Heisenberg models, where the state
space is the $d$--dimensional unit sphere  
and the interactions are proportional to
the cosines of the angles between neighboring spins.  The phenomenon of
interest here is the classification of phase transition
(non-uniqueness of the Gibbs state) according to whether it is {\it robust}.
In many cases, including all of the Heisenberg and Potts models,
occurrence of robust phase transition is determined by the geometry 
(branching number) of
the tree in a way that parallels the situation with independent 
percolation and usual phase transition for the Ising model. The critical
values for robust phase transition for the Heisenberg and Potts
models are also calculated exactly. In some cases, such as the $q\ge 3$ Potts 
model, robust phase transition and usual phase transition do not coincide, 
while in other cases, such as the Heisenberg models, we conjecture that robust
phase transition and usual phase transition are equivalent. In addition,
we show that symmetry breaking is equivalent to the existence of a phase
transition, a fact believed but not known for the rotor model on $\Z^2$.
\end{abstract}

Running head: Phase transitions for Heisenberg models.
\end{titlepage}

\setcounter{equation}{0}

\section{Definition of the model and main results} \label{sec:one}

Particle systems on trees have produced the first and most tractable 
examples of certain qualitative phenomena.  For example, the
contact process on a tree has multiple phase transitions,
(\cite{Pem,Lig2,Sta}) and the critical temperature for the
Ising model on a tree is determined by its branching number or Hausdorff 
dimension (\cite{Ly1,EKPS,PP}), which makes the Ising model 
intimately related to independent percolation whose critical value is also
determined by the branching number (see \cite{Ly2}).  In this paper
we study several models on general infinite trees,
including the classical Heisenberg and Potts models.
Our aim is to exhibit a distinction between two kinds of phase transitions,
{\it robust} and {\it non-robust}, as well as to investigate conditions
under which robust phase transitions occur.  

In many cases, including the Heisenberg and Potts models,
the existence of a robust phase transition is determined 
by the branching number.  However, in some cases (including
the $q > 2$ Potts model), the critical temperature for 
the existence of usual phase transition is not determined by the 
branching number.  Thus robust phase transition 
behaves in a more universal manner than non-robust phase transition, 
being a function of the branching number alone,
as it is for usual phase transition for independent percolation and
the Ising model.  Although particle systems on trees do not always
predict the qualitative behavior of the same particle system on 
high-dimensional lattices, it seems likely that there is a
lattice analogue of non-robust phase transition, which would make
an interesting topic for further research.  Another unresolved question
is whether there is ever a non-robust phase transition for the Heisenberg
models (see Conjecture~\ref{conj:PS}).  

We proceed to define the general statistical ensemble 
on a tree and to state the main results of the paper.  
Let $G$ be a compact metrizable group acting transitively by isometries on 
a compact metric space $(\S,d)$.
It is well known that there exists a unique
$G$--invariant probability measure on $\S$, which we denote by $\dx$.  
An {\bf energy function} is any nonconstant function 
$H : \S \times \S \rightarrow \R$ that is symmetric, continuous,
and $d$--invariant in that $H(x,y)$ depends only on $d(x,y)$. This implies that
$$
H(x,y)=H(gx,gy) \,\,\forall \, x,y\in \S, \,g\in G.
$$

$\S$ together with its $G$--action and the function $H$ will be called a
{\bf statistical ensemble}.
Several examples with which we will be concerned are as follows.

\begin{eg} \label{eg:ising}
The Ising model.  Here $\S = \{ 1 , -1 \}$
acted on by itself (multiplicatively), $d$ is the usual discrete metric,
$\dx$ is uniform on $\S$, and $H (x , y) = - xy$.  
\end{eg}

\begin{eg} \label{eg:potts}  
The Potts model.  Here $\S = \{ 0 , 1 , \ldots ,
q-1 \}$ for some integer $q > 1$, $G$ is the symmetric group $S_q$ with its
natural action, $d$ is the usual discrete metric,
$\dx$ is uniform on $\S$, and $H (x , y) = 1 - 2 \delta_{x,y}$.
This reduces to the Ising model when $q = 2$.
\end{eg}

\begin{eg} \label{eg:rotor}
The rotor model.  Here $\S$ is the unit 
circle, acted on by itself by translations, 
$d(\theta , \phi) = 1- \cos (\theta - \phi)$,
$\dx$ is normalized Lebesgue measure, and 
$H (\theta , \phi) = - \cos (\theta - \phi)$.  
\end{eg}

\begin{eg} \label{eg:spherical} 
The Heisenberg models for $d \ge 1$.
In the $d$--dimensional
Heisenberg model, $\S$ is the unit sphere $S^d$, $G$ is
the special orthogonal group with its natural action; $d(x,y)$ is $1-x\cdot y$,
$\dx$ is normalized surface measure, 
and $H (x , y)$ is again the negative of the dot product of $x$ and $y$. 
When $d=1$, we recover the rotor model.
\end{eg}

Let $A$ be any finite graph, with vertex and edge sets denoted by
$V(A)$ and $E(A)$ respectively, and let $\J : E (A) \rightarrow \R^+$
be a function mapping the edge set of $A$ to the nonnegative reals which
we call {\bf interaction strengths}.
We now assume that $\S$, $G$ and $H$ are given and fixed.

\begin{defn} \label{defn:Gibbs}
The {\bf Gibbs measure} with interaction strengths $\J$ is
the probability measure $\mu=\mu^{\J}$ on $\S^{V(A)}$ whose density
with respect to product measure $\dx^{V(A)}$ is given by
$$ { \exp (- H^\J (\eta)) \over Z},\,\,\,\, \eta\in \S^{V(A)} $$
where
$$ H^\J(\eta) = \sum_{e = \overline{xy} \in E(A)} 
   \J (e) H(\eta (x) , \eta (y)) ,$$
and $Z = \int \exp (- H^\J(\eta)) \, \dx^{V(A)}$ is a normalization. 
\end{defn}

In statistical mechanics, one wants to define Gibbs measures on infinite
graphs $A$ in which case the above definition of course does not make sense.
We follow the usual approach (see~\cite{Ge}), in which one introduces
boundary conditions and takes a weak limit of finite subgraphs
increasing to $A$.  Since the precise nature of the boundary conditions 
play a role here (we know this to be true at least for the Potts model 
with $q > 2$), we handle boundary conditions with extra care and, 
unfortunately, notation.  We give definitions in the case of a rooted 
tree, though the extensions to general locally finite graphs are immediate.  
By a {\bf tree}, we mean any connected loopless graph $\Gamma$ where
every vertex has finite degree.  One fixes a vertex $o$ of $\Gamma$ 
which we call the {\bf root}, obtaining a {\bf rooted tree}.  
The vertex set of $\Gamma$ is denoted by $V(\Gamma)$.  If $x$ is a vertex,
we write $|x|$ for the number of edges on the shortest
path from $o$ to $x$ and for two vertices $x$ and $y$, we
write $|x-y|$ for the number of edges on the shortest
path from $x$ to $y$.  For vertices $x$ and $y$, we write $x \le y$ 
if $x$ is on the shortest path from $o$ to $y$, $x < y$
if $x \le y$ and $x \ne y$, and $x \to y$ if $x \le y$ and 
$|y|=|x|+1$.  For $x \in V(\Gamma)$, the tree $\Gamma (x)$ denotes
the subtree of $\Gamma$ rooted at $x$ consisting of $x$ and all of its
descendents. We also define $\partial\Gamma$, which we refer to as the
boundary of $\Gamma$, to be the set of infinite
self-avoiding paths starting from $o$.
Throughout the paper, the following assumption is in force.

\noindent{\bf ASSUMPTION:} For all trees considered in this paper,
the number of children of the vertices will be assumed bounded
and we will denote this bound by $B$.

A {\bf cutset} $C$ is a finite set of vertices not including $o$ such that
every self-avoiding infinite path from $o$ intersects $C$ and such that there
is no pair $x , y \in C$ with $x < y$. Given a cutset $C$,
$\Gamma \backslash C$ has one finite component (which contains $o$)
which we denote by $C^i$ (``i'' for inside) and we let 
$C^o$ (``o'' for outside) denote the union of the infinite components
of $\Gamma \backslash C$. We say that a sequence
$\{ C_n \}$ of cutsets approaches $\infty$ if for all $v \in \Gamma$,
$v \in C_n^i$ for all sufficiently large $n$.  

Boundary conditions will take the form of specifications of the
value of $\eta$ at some cutset $C$. 
Let $\delta$ be any element of $\S^C$.  The Gibbs measure with boundary
condition $\delta$ is the probability 
measure $\mu^\delta_C = \mu^{\J , \delta}_C$ on $\S^{C^i}$
whose density with respect to product measure $\dx^{C^i}$ is given by
\begin{equation} \label{eq:Gibbs}
{ \exp (- H^{\J , \delta}_C (\eta)) \over Z},\,\,\,\, \eta\in \S^{C^i} 
\end{equation}
where
$$ H^{\J , \delta}_C (\eta) = \sum_{e = \overline{xy} \in E(\Gamma) 
   \atop x,y \in C^i} \J (e) H(\eta (x) , \eta (y)) 
   + \sum_{e = \overline{xy} \in E(\Gamma) \atop x \in C^i, y \in C} 
   \J (e) H(\eta (x) , \delta (y)) $$
and 
$Z = \int \exp (- H^{\J , \delta}_C (\eta)) \, \dx^{C^i}$ is a normalization. 
When we don't include the second summand above, 
we call this the {\it free} Gibbs measure on $C^i$, denoted by $\mufree_C$,
where $\J$ is suppressed in the notation.
As we will see in Lemma~\ref{lem:free}, the free measure 
does not depend on $C$ except for its domain of definition, 
so we can later also suppress $C$ in the notation.

\begin{defn} \label{defn:gibbstree}
A probability measure $\mu$ on $\S^{V(\Gamma)}$ is called a {\bf Gibbs state} 
for the interactions $\J$ if for each cutset $C$, the 
conditional distribution on  $C^i$ given the configuration $\delta'$ 
on $C\cup C^o$  is given by $\mu_C^{\J , \delta}$ where $\delta$ 
is the restriction of $\delta'$ to $C$. (A similar definition is used for
general graphs.)
Both in the case of lattices and trees (or for any graph), we say that
a statistical ensemble {\bf exhibits a phase transition (PT) for the
interaction strengths $\J$} if there is more than one Gibbs state
for the interaction strengths $\J$.
\end{defn}

In the next section we will prove
\begin{lem} \label{lem:free}
Fix interaction strengths $\J$ and let 
$C$ and $D$ be any two cutsets of $\Gamma$.  Then the projections
of $\mufree_C$ and $\mufree_D$ to $\S^{C^i \cap D^i}$ are equal. Hence the 
measures $\mufree_C$ have a weak limit as $C \rightarrow \infty$,
denoted $\mufree$.
\end{lem}

For general graphs, the measures $\mufree_C$ are not compatible 
in this way. Also, one has the following fact, which follows from 
Theorems~4.17 and~7.12 in \cite{Ge}.

\begin{lem} \label{lem:limits}
If $\{C_n\}$ is a sequence of cutsets approaching $\infty$ and if for each
$n$, $\delta_n\in \S^{C_n}$, then any weak subsequential limit of the sequence 
$\{\mu_{C_n}^{\J,\delta_n}\}_{n\ge 1}$
is a Gibbs state for the interactions $\J$. 
In addition, if all such possible limits are the same,
then there is no phase transition.  (A similar statement holds for 
graphs other than trees.)
\end{lem}

We pause for a few remarks about more general graphs, before restricting
our discussion to trees for the rest of the paper.  Lemma~\ref{lem:free}
does not apply to graphs with cycles, so the existence of a unique
weak limit $\mufree$ is not guaranteed there, but Lemma~\ref{lem:limits}
together with compactness tells us that there always is at least one Gibbs
state.  The state of knowledge about the rotor model (Example~\ref{eg:rotor}) 
on more general graphs is somewhat interesting.  It is known 
(see~\cite{Ge}, p.178 and p.434)
that for $\Z^d$, $d \leq 2$, all Gibbs states are rotationally invariant
when $\J\equiv J$ for any $J$ 
(and it is believed but not known that there is a unique Gibbs state for the 
rotor model in this case) while for $d \geq 3$, there are values of $J$ 
for which the rotor model with $\J\equiv J$ has a Gibbs state whose
distribution at the origin is not
rotationally invariant (and hence there is 
more than one Gibbs state). In statistical mechanics,
this latter phenomenon is referred to as a {\it continuous 
symmetry breaking} since we have a continuous state space (the circle) where
the interactions are invariant under a certain continuous symmetry 
(rotations) but there are Gibbs states which are not invariant under this 
symmetry.  We also mention that it is proved in~\cite{C} that 
for the rotor model with $\J\equiv J$ for any $J$ 
on any graph of bounded degree  
for which simple random walk is recurrent, all
the Gibbs states are rotationally invariant. (This was then extended
in~\cite{MW} where the condition of boundedness of the degree is dropped
and the group involved is allowed to be more general than the circle.)
This however is not a sharp criterion: in~\cite{E}, 
a graph (in fact a tree) is constructed for which simple random walk 
is transient but such that there is no phase transition in the 
rotor model when $\J\equiv J$ for any $J$. (This will also follow from
Theorem~\ref{th:0hd} below together with the easy fact that
there are trees with branching number 1 for which
simple random walk is transient.)
However, Y.\ Peres has conjectured a sharp criterion,
Conjecture~\ref{conj:peres} below, for which our Corollary~\ref{cor:perestree} 
together with the discussion following it provides some corroboration. 

For the rest of this paper, we will restrict to trees.  It is usually 
in this context that the most explicit results can be obtained and our 
basic goal is to determine whether there is a phase transition by comparing 
the interaction strengths with the ``size'' (branching number) of our tree. 
It turns out that we can only 
partially answer this question but the question which we can answer more
completely is whether there is a {\it robust} phase transition, a concept 
which we will introduce shortly.

\begin{defn} \label{defn:notation}
Given $\J,C$ and $\delta$ defined on $C$, let
$f^{\J,\delta}_{C , o}$ 
(or $f^{\delta}_{C , o}$ if $\J$ is understood)
denote the marginal density of $\mu^{\J , \delta}_{C}$ at the root $o$.  
\end{defn}

For any tree, recall that $\Gamma (v)$ denotes the subtree 
rooted at $v$, so that the tree $\Gamma(v)$ has vertex set 
$\{ w \in \Gamma : v \leq w\}$. If $v\in C^i$ and
we intersect $C$ with $\Gamma (v)$, we obtain a cutset $C(v)$ for $\Gamma(v)$.
We now extend Definition~\ref{defn:notation} to other marginals as follows.

\begin{defn} \label{defn:marginals}
With $\J,C$ and $\delta$ as in Definition~\ref{defn:notation} and $v \in C^i$, 
define $f_{C,v}^{\J , \dd}$ 
by replacing $\Gamma$ by $\Gamma (v)$, $C$ with $C(v)$, $\J$ with
$\J$ restricted to $E (\Gamma (v))$, $\delta$ with $\delta$ restricted
to $C(v)$ and $o$ with $v$ in Definition~\ref{defn:notation}.  
\end{defn}

It is important to note that $f_{C,v}^{\J,\dd}$ is not the density 
of the projection of $\mu_C^{\J,\dd}$ onto vertex $v$, but rather
the density of a Gibbs measure with similar boundary conditions
on the smaller graph $\Gamma (v)$.

\begin{defn} \label{defn:SB}
A statistical ensemble on a tree $\Gamma$ exhibits a 
{\bf symmetry breaking (SB) for the interactions $\J$} 
if there exists a Gibbs state such that the marginal distribution at some
vertex $v$ is not $G$--invariant (or equivalently is not $\dx$).
\end{defn}

The following proposition which will be proved in Section~\ref{sec:prelims} 
is interesting since it establishes the equivalence of PT and SB for
general trees and general statistical ensembles, something not known for 
general graphs, see the remark below.

\begin{prop} \label{prop:SB=}
Consider a statistical ensemble on a tree $\Gamma$ with interactions $\J$. 
The following four conditions are equivalent. \\
(i) There exists a vertex $v$ such that for 
any sequence of cutsets $C_n\to\infty$, there exist boundary 
conditions $\p_n$ on $C_n$ such that
$$
\inf_n \|f_{C_n,v}^{\p_n}-1\|_\infty \neq 0.
$$
(ii)  There exists a vertex $v$, a
sequence of cutsets $C_n\to\infty$ and boundary 
conditions $\p_n$ on $C_n$ such that
$$
\inf_n \|f_{C_n,v}^{\p_n}-1\|_\infty \neq 0.
$$
(iii) The system satisfies SB. \\
(iv) The system satisfies PT.
\end{prop}

We now fix a distinguished element in $\S$, hereafter denoted $\zero$.  
The notation $\mu^{\J , +}_C$ denotes $\mu^{\J , \delta}_C$ when 
$\delta$ is the constant function $\zero$.  In the case $\J \equiv J$, we
denote this simply $\mu^{J,+}_C$.  We will be particularly concerned
about whether $\mu^{\J , +}_C \rightarrow \mufree$ weakly, as 
$C \rightarrow \infty$.    

\begin{defn} \label{defn:SB+}
A statistical ensemble on a tree $\Gamma$ exhibits a 
{\bf symmetry breaking with plus boundary conditions (SB+) for the 
interactions $\J$} if there exists a vertex $v$ and a sequence of 
cutsets $C_n\to\infty$ such that
$$
\inf_n \|f_{C_n,v}^{\J,+}-1\|_\infty \neq 0.
$$
\end{defn}

Note that by symmetry, SB+ does not depend on which point of $\S$
is chosen to be $\zero$.

In Section~\ref{sub:spherical} we will prove:
\begin{pr} \label{pr:rotor equiv}
For the rotor model on a tree, SB is equivalent to SB+.  
\end{pr}

We conjecture but cannot prove the stronger statement:
\begin{conj} \label{conj:SB}
For any Heisenberg model on any graph, SB is equivalent to SB+.
\end{conj}

\noindent{\em Remarks:} $(i)$ By Proposition~\ref{prop:SB=}, we have that 
SB+ implies SB for any statistical ensemble on a tree.  While 
Proposition~\ref{prop:SB=} tells us that PT and SB are equivalent for any
statistical ensemble
on a tree, we note that such a result is not even known for the rotor model
on $\bbbz^2$ where it has been established that for all
$J$, all Gibbs states are rotationally invariant for $\J\equiv J$
but where it has not been established that there is no phase transition. 
A weaker form of the above conjecture would be that SB+ and SB are 
equivalent for all Heisenberg models on trees.  This is
Problem~\ref{pblm:all spheres} in~Section~\ref{sec:anal}.
An extension to graphs with cycles would seem to entail a different kind
of reasoning, perhaps similar to the inequalities of Monroe and 
Pearce~\cite{MP} which fall just short of proving
Conjecture~\ref{conj:SB} for the rotor model. \\
\noindent{$(ii)$} The fact that PT and SB+ 
are equivalent when the rotor model is replaced by the Ising model 
is an immediate consequence of the fact that the probability measure is 
stochastically increasing in the boundary conditions. More generally,
it is also the case that PT and SB+ are equivalent for the Potts models
(see \cite{ACCN}).

We now consider the idea of a {\it robust phase transition} where we 
investigate if the boundary conditions on a cutset have a nontrivial effect
on the root even when the interactions along the cutset are made
arbitrarily small but fixed.

Given parameters $J>0$ and $J' \in (0,J]$ and
a cutset $C$ of $\Gamma$, let $ \J ( J', J , C)$ be the function 
on $E(\Gamma)$ which is $J$ on edges in $C^i$ and
$J'$ on edges connecting $C^i$ to $C$ (the values elsewhere being
irrelevant).  
Let $f^{J',J , +}_{C_n , o}$ denote the marginal at the root $o$ of the
measure $\mu^{J',J, +}_C:=\mu^{\J (J', J , C) , +}_C$.

\begin{defn} \label{defn:robustPT}
The statistical ensemble on the tree $\Gamma$ has a {\bf robust phase 
transition (RPT) for the parameter $J>0$} if for every $J'\in (0,J]$
$$
\inf_C \|f^{J',J , +}_{C , o} - 1\|_\infty \neq 0 \, 
$$
where the $\inf$ is taken over all cutsets $C$.
\end{defn}

\noindent{\em Remarks:}  
In the case $\J \equiv J$, by taking $J'=J$, it is clear
that a RPT implies SB+ (which in turn implies SB and PT).  Note that
in this case, RPT is
stronger than SB+ not only because $J'$ can be any number in $(0,J]$
and the root $o$ must play the role of $v$
but also because in SB+, we only require that for {\it some} sequence
of cutsets going to infinity, the marginal at the vertex $v$ stays away from
uniform while in RPT, we require this for {\it all} cutsets going to
infinity.
We note also that with some care, this definition makes sense for general 
graphs, and that the issue of robustness of phase transition on general graphs
is worth investigating, although we do not do so here.

Our first theorem gives criteria based on $J$ and the branching number
of $\Gamma$ (which will now be defined) for robust phase transition to
occur for the Heisenberg models. A little later on, we will have an analogous
result for the Potts models. In \cite{F}, Furstenberg introduced the notion
of the Hausdorff dimension of a tree (or more accurately of the boundary of 
the tree). This was further investigated by Lyons~(\cite{Ly2}) using the 
term branching number instead.
The {\bf branching number} of a tree $\Gamma$, denoted
$\textstyle {br} (\Gamma)$, is a real number greater than or equal to one
that measures the average number of branches per vertex of the
tree. More precisely, the {\bf branching number} of $\Gamma$ is defined by
$$\textstyle{br}\,\Gamma:=\inf\left\{\lambda>0;\inf\limits_{C}
\sum_{x \in C}\lambda^{-|x|} = 0 \right \} \;$$
where the second infimum is over all cutsets $C$.
The branching number is a measure of the average number of
branches per vertex of $\Gamma$.
It is  less than or equal to $\liminf_{n \to \infty} M_n^{1/n}$, where
$M_n := | \left \{x \in \Gamma ; |x| = n \right \}|$,
and takes more of the structure of $\Gamma$ into account
than does this latter growth rate.
For sufficiently regular trees, such as homogeneous trees or, more generally,
Galton-Watson trees, $\textstyle{br}\, \Gamma = \lim_{n\to\infty}
M_n^{1/n}$ (\cite{Ly2}). We also mention that the branching number is the
exponential of the Hausdorff dimension of $\partial\Gamma$ where the latter
is endowed with the
metric which gives distance $e^{-k}$ to two paths which split off after $k$
steps. As indicated earlier, the branching number has been an important
quantity in previous investigations. More specifically, in \cite{Ly1} and
\cite{Ly2}, the critical values for independent percolation and
for phase transition in the Ising model on general trees 
are explicitly computed in terms of the branching number.

For each $J\ge 0$, define a continuous strictly positive 
probability density function 
$K_J : \S \rightarrow \R^+$ by 
\begin{equation} \label{eq:KJ}
K_J (u): = C(J)^{-1} \exp (- J H (u , \zero)) 
\end{equation}
where $C(J) = \int \exp (- J H(w,\zero)) \, \dx (w)$ is a normalizing
constant, and more generally let
$K_{J,y} : \S \rightarrow \R^+$ be given by 
\begin{equation} \label{eq:KJy}
K_{J,y} (u): = C(J)^{-1} \exp (- J H (u , y)) 
\end{equation}
(noting that $K_{J,\zero}=K_{J}$).
Let $\KJ$ denote the convolution operator on the space
$L^2 (\S , \dx)$ given by the formula
\begin{equation} \label{eq:conv}
 \KJ f (u) : = \int_{\S} f(x) K_{J,x} (u) \dx(x) \,\, .
\end{equation}
Note that by the assumed invariance $\int_{\S} \exp (- J H(w,y)) \, \dx (w)$
is independent of $y$ and that $f\ge 0$ and $\int_{\S} f(x) \dx(x)=1$ imply
that $\KJ f\ge 0$ and $\int_{\S} \KJ f(x) \dx(x)=1$.
We extend the above notation to cover the case where $f$ is a pointmass $\p_y$ 
at $y$ by defining in that case 
\begin{equation} \label{eq:pointconv}
\KJ \p_y (u) : = K_{J , y}(u) .
\end{equation}

%%\noindent{\em Remark:}  In particular, we see for the 
%%spherical models (and will see for the Potts  models) that
%%the critical value for $J$ for robust phase transition is determined by 
%%$\textstyle {br}  (\Gamma)$ (equivalently, 
%%there is a critical branching number for each $J$ (and $d$)).

We will now give the exact critical
parameter $J$ for RPT for the Heisenberg models. For any $d\ge 1$, let 
$$
\calrdJ:= {\int_{-1}^1 r e^{Jr}(1-r^2)^{{d \over 2}-1} dr \over
\int_{-1}^1  e^{Jr}(1-r^2)^{{d \over 2}-1} dr }.
$$

When $d=1$ (rotor model), this is (by a change of variables)
the first Fourier coefficient  of $K_J$
($\int_{\S} K_J(\theta) \cos (\theta) d\theta$) which is perhaps
more illustrative.
When $d=2$, this is the first Legendre coefficient 
of $e^{Jr}$ (properly normalized) and for $d\ge 3$,
this is the first so-called ultraspherical coefficient of $e^{Jr}$ 
(properly normalized).

\begin{th} \label{th:main} Let $d\ge 1$. \\
(i) If $\textstyle {br} (\Gamma) \calrdJ <1$,
then the $d$--dimensional Heisenberg model on $\Gamma$ with
parameter $J$ does not exhibit a robust phase transition. \\
(ii) If $\textstyle{br}(\Gamma) \calrdJ >1$, then the $d$--dimensional
Heisenberg model on $\Gamma$ with parameter $J$ exhibits a robust phase 
transition.
\end{th}

\noindent{\em Remark:}
It is easy to see that $\lim_{d\to\infty} \calrdJ=0$ which says that 
it is harder to obtain a robust phase transition on higher dimensional
spheres. This is consistent with the fact that it is in some sense harder to
have a phase transition for the rotor model than in the Ising 
model (0-dimensional sphere); this
latter fact can be established using the ideas in \cite{PS}.

A simple computation shows that the derivative of $\calrdJ$ with respect to
$J$ is the variance of a random variable whose density function is 
proportional to $e^{Jr}(1-r^2)^{d/2 -1}$ on $[-1,1]$, thereby obtaining
the following lemma.
\begin{lem} \label{lem:inc}
For any $d\ge 1$,
we have that $\calrdJ$ is a strictly increasing function of $J$.
\end{lem}

Theorem \ref{th:main} and Lemma \ref{lem:inc} together with the
fact that for any $d\ge 1$,
$\calrdJ$ is a continuous function of $J$ which approaches
0 as $J\to 0$ and approaches 1 as $J\to \infty$ give us
the following corollary.

\begin{cor} \label{cor:critical}
For any Heisenberg model with $d\ge 1$ and any
tree $\Gamma$ with branching number larger than 1, let $J_c=J_c(\Gamma, d)$
be such that $\textstyle{br}(\Gamma) \rho^d(J_c)=1$. Then
there is a robust phase transition for the $d$--dimensional Heisenberg 
model on $\Gamma$ if $J> J_c$ and
there is no such robust phase transition for $J< J_c$.
\end{cor}

For the Heisenberg models, we believe that phase transition and
robust phase transition coincide and therefore we have the
following conjecture.

\begin{conj} \label{conj:PS}
For any $d\ge 1$, if $\textstyle {br} (\Gamma) \calrdJ <1$, 
then the $d$--dimensional Heisenberg 
model on $\Gamma$ with parameter $J$ does not exhibit a phase 
transition.
\end{conj}

We can however obtain the following weaker form of this conjecture which is
valid for all statistical ensembles.

\begin{th} \label{th:0hd}
If $\textstyle {br}  (\Gamma) = 1$, then there is no phase transition
for any statistical ensemble on $\Gamma$ with bounded $\J$.
\end{th}

Theorems \ref{th:main}(ii) and \ref{th:0hd} together with the facts that
RPT implies PT and that for any $d\ge 1$,
$\lim_{J\to\infty} \calrdJ = 1$ immediately yield the following corollary.

\begin{cor} \label{cor:perestree}
For any Heisenberg model with $d\ge 1$ and for any tree $\Gamma$, 
there is a \pt for the tree $\Gamma$ for some
value of the parameter $J$ if and only if $\textstyle {br} (\Gamma)>1$.
\end{cor}

Since it is known (see \cite{Ly2}) that $\textstyle {br} (\Gamma)>1$
if and only if there is some $p< 1$ with the property that when
performing independent percolation on $\Gamma$ with parameter $p$, there
exists a.s.\ an infinite cluster on which simple random walk is transient, the
above corollary yields the following conjecture of Y. Peres for the
special case of trees of bounded degree.

\begin{conj} \label{conj:peres}
For any graph $A$, the rotor model exhibits a \pt for some $J$ if and only
if there is some $p< 1$ with the property
that performing independent bond percolation on $A$ with parameter $p$, there
exists a.s.\ an infinite cluster on which simple random walk is transient.
\end{conj}

Recall that the rotor model on the graph $A$ exhibits
no SB for any parameter $J$ if $A$ is recurrent for simple random walk,
which is of course consistent with the above conjecture.
Note that, on the other hand, the standard Ising model does exhibit a
\pt on $\Z^2$, a graph which is recurrent (as are its subgraphs)
for simple random walk.

The next result states the critical value for RPT for the Potts models.

\begin{th} \label{th:potts}
Consider the Potts model with $q \ge 2$ and let
$$\alpha_J = {e^J - e^{-J} \over e^J + (q-1) e^{-J}} \, .$$
(i) If $\textstyle {br} (\Gamma) \alpha_J <1$,
then the Potts model on $\Gamma$ with
parameter $J$ does not exhibit a robust phase transition. \\
(ii) If $\textstyle{br}(\Gamma) \alpha_J >1$, then the 
Potts model on $\Gamma$ with parameter $J$ exhibits a robust phase 
transition.
\end{th}

\noindent{\em Remarks:}\\
$(i)$ $d\alpha_J/dJ >0$ and so there is a critical value of $J$ depending
on $\textstyle{br}(\Gamma)$ analogous to in Corollary~\ref{cor:critical}
for the Heisenberg models. \\
\noindent{$(ii)$} 
Note that when $q=2$ (the Ising model), this formula agrees 
with the formula for the Heisenberg models when one formally sets $d=0$
in the formula 
$$
\calrdJ=\int_{S^d} (x \cdot \zero) K_J (x) \, \dx(x),
$$
the latter being obtained by a change of variables.

To point out the subtlety involved in Conjecture \ref{conj:PS}, 
we continue to discuss the Potts model, a case in which the analogue of 
Conjecture~\ref{conj:PS} fails. Our final result tells us that 
phase transitions (unlike robust phase transitions) in the Potts model with
$q > 2$ cannot be determined by the branching number.

\begin{th} \label{th:2trees} Given any integer $q >2$,
there exist trees $\Gamma_1$ and $\Gamma_2$ and a nontrivial interval $I$ such
that $\textstyle {br} (\Gamma_1) < \textstyle {br} (\Gamma_2)$ 
and for any $J\in I$, there is a phase transition for the $q$--state
Potts model with parameter $J$ on $\Gamma_1$ but 
no such phase transition on $\Gamma_2$.
\end{th}

\newpage
\noindent{\em Remarks:}\\
$(i)$ $\Gamma_1$ and $\Gamma_2$ can each be taken to
be spherically symmetric which means that for all $k$, all vertices at the
$k$th generation have the same number of children. \\
\noindent{$(ii)$} 
In the case $q=2$, more is known.  In \cite{Ly1}, the critical value
for phase transition in the Ising model is found and corresponds to
what is obtained in Theorem~\ref{th:potts} above. It follows that 
there is never
a non-robust phase transition except possibly at the critical value. However,
a sharp capacity criterion exists~\cite{PP} for phase transition for the
Ising model (settling the issue of phase transition at the critical parameter)
and using this criterion, one can show that phase transition and robust 
phase transition correspond even at criticality.
The arguments of~\cite{PP} cannot be extended to the Potts model for $q > 2$
because the operator $\KJ$, acting on a certain likelihood function,
when conjugated by the logarithm is not concave in this case.
Theorems~\ref{th:potts} and~\ref{th:2trees}
together tell us that there is indeed a non-robust phase
transition when $q > 2$ for a nontrivial interval of $J$.

The rest of the paper is devoted to the proofs of the above results.
In Section~\ref{sec:prelims}, we collect several lemmas that apply to
general statistical ensembles, including the basic recursion 
formula (Lemma~\ref{lem:rec}) that allows us to analyze general statistical 
ensembles on trees, prove Lemma~\ref{lem:free} and
Proposition~\ref{prop:SB=} as well as provide some background
concerning Heisenberg models (showing that they satisfy the more
general hypotheses of Theorems~\ref{th:gen ii} and~\ref{th:gen i} given
later on) and the more general notion of distance
regular spaces.  Section~\ref{sec:proofs} is devoted to the proofs of
Theorems~\ref{th:gen ii} and~\ref{th:gen i}.  In
Section~\ref{sec:anal}, we use these theorems to find the critical
parameters for robust phase transition in the Heisenberg and
Potts models, Theorems~\ref{th:main} and~\ref{th:potts}, as well as
prove Proposition~\ref{pr:rotor equiv}.
Section~\ref{sec:zero} discusses the special case of trees of
branching number 1, proving Theorem~\ref{th:0hd}. Finally, in
Section~\ref{sec:potts}, Theorem~\ref{th:2trees} is proved.

\setcounter{equation}{0}

\section{Basic background results} \label{sec:prelims}
In this section, we collect various background results which
will be needed to prove the results described in the introduction.
We begin with a subsection describing results pertaining to trees that 
hold for general statistical ensembles. 
After discussing the concept of a distance
regular space in Section~\ref{sub:drs}, we specialize 
to Heisenberg models (the most relevant family of continuous distance 
regular models) in Section~\ref{sub:sphere}
and then to distance regular graphs in Section~\ref{sub:finite}.  

\subsection{The fundamental recursion and other lemmas} \label{sub:rec}

We start off with two lemmas exploiting the recursive structure
of trees.

Let $\S,G$ and $H$ be a statistical ensemble.
Let $A_1$ and $A_2$ be two disjoint finite graphs, with distinguished 
vertices $v_1 \in V(A_1)$ and $v_2 \in V(A_2)$.  Let $\J_1$ and $\J_2$ 
be interaction functions for $A_1$ and $A_2$,
i.e., positive functions on $E(A_1)$ 
and $E(A_2)$ respectively.  For any $C_1 \subseteq V(A_1) \setminus 
\{ v_1 \}$ (possibly empty) and any $C_2 \subseteq V(A_2)$, and for any 
$\delta_1 \in \S^{C_1}$ and $\delta_2 \in \S^{C_2}$, 
we have measures $\mu_i := \mu^{\J_i , \delta_i}_{C_i}$, $i = 1, 2$
on $\S^{V(A_i)\setminus C_i}$ defined (essentially) by~(\ref{eq:Gibbs}). 
Abbreviate $H_{C_i}^{\J_i,\p_i}$ (which has the obvious meaning) by $H_i$.
Let $A$ be the union of $A_1$
and $A_2$ together with an edge connecting $v_1$ and $v_2$.  Let
$C = C_1 \cup C_2$, $\J$ extend each $\J_i$ and the value of the
new edge be given the value $J$, $\delta$
extend each $\delta_i$ and denote $\mu^{\J , \delta}_C$ 
(a probability measure on
$\S^{(V(A_1)\setminus C_1)\cup (V(A_2)\setminus C_2)}$)
by $\mu$ and $H^{\J , \delta}_C$ (again having the obvious meaning) by $H$.
The identity \begin{equation} \label{eq:Hdecomp}
H = H_1 + H_2 + J H (\eta (v_1) , \eta (v_2))
\end{equation}
leads to the following lemma.

\begin{lem} \label{lem:decomps}
The measure $\mu$ satisfies
\begin{equation} \label{eq:mudecomp}
{d\mu \over d (\mu_1 \times \mu_2)} = c \exp [- J H(\eta_1 (v_1) , 
   \eta_2 (v_2))] ,
\end{equation}
where 
$$c = \left [ \int \int \exp (- J H (\eta_1 (v_1) , \eta_2 (v_2))) 
   \, d\mu_1 (\eta_1) \, d\mu_2 (\eta_2) \right ]^{-1}$$
is a normalizing constant.  Let $f_i$ denote the marginal density
of $\mu_i$ at $v_i$, $i = 1 , 2$, and $f$ denotes the marginal
density of $\mu$ at $v_1$.  Then the projection $\mu^{(1)}$ of
$\mu$ onto $\S^{V(A_1)\setminus C_1}$ satisfies
\begin{equation} \label{eq:mudecomp2}
\mu^{(1)} = c \int \int \mu_{1 , y} f_1 (y) f_2 (z) \exp (- J H(y,z))
   \, \dx (z) \, \dx (y)
\end{equation}
for some normalizing constant $c$, where $\mu_{1 , y}$ denotes the
conditional distribution of $\mu_1$ given $\eta (v_1) = y$.  Consequently,
\begin{equation} \label{eq:fdecomp}
f (y) = c f_1 (y) \int f_2 (z) \exp (- J H(y,z)) \, \dx (z) \, ,
\end{equation}
where $c$ normalizes $f$ to be a probability density.  
\end{lem}

\noindent{\bf Proof.}  The relation~(\ref{eq:mudecomp}) follows
from~(\ref{eq:Hdecomp}) and the defining equation~(\ref{eq:Gibbs}).  From 
this it follows that the measure $\mu$ on pairs $(\eta_1 , \eta_2)$
makes $\eta_1$ and $\eta_2$ conditionally independent given $\eta_1 (v_1)$
and $\eta_2 (v_2)$.  Hence the conditional distribution of $\mu^{(1)}$
given $\eta_1 (v_1) = y$ and $\eta_2 (v_2) = z$ is just $\mu_{1 , y}$.
Next, (\ref{eq:mudecomp}) and the last fact
yield~(\ref{eq:mudecomp2}).  The marginal of 
$\mu_{1, y}$ at $v_1$ is just $\delta_y$, and
so~(\ref{eq:mudecomp2}) yields~(\ref{eq:fdecomp}).  
$\qed$

A tree $\Gamma$ may be built up from isolated vertices by the joining 
operation described in the previous lemma.  The decompositions in 
Lemma~\ref{lem:decomps} may be applied inductively 
to derive a fundamental recursion for marginals.  This recursion, 
Lemma~\ref{lem:rec} below, expresses the marginal distribution
at the root of $\Gamma$ as a pointwise product of marginals at
the roots of each of the generation 1 subtrees, each convolved
with a kernel $K_J$.  The normalized pointwise product will be
ubiquitous throughout what follows, so we introduce notation for it.

\begin{defn}
If $f_1 , \ldots , f_k$ are nonnegative functions on $\S$ with
$\int f_i \, \dx = 1$ for each $i$, 
let \hfill \\ $\poi_k (f_1 , \ldots , f_k)$ 
denote the normalized pointwise product, 
$$\poi_k (f_1 , \ldots , f_k) (x) = {\prod_{i=1}^k f_i (x) \over
   \int \prod_{i=1}^k f_i (y) \, \dx (y)}$$
whenever this makes sense, e.g., when each $f_i$ is in $L^k (\dx)$
and the product is not almost everywhere zero.  Let $\poi$ denote
the operator which for each $k$ is $\poi_k$ on each $k$-tuple of functions.
There is an obvious associativity property, namely $\poi (\poi (f,g) , h)
= \poi (f,g,h)$, which may be extended to arbitrarily many arguments.
\end{defn}

\begin{lem}[Fundamental recursion] \label{lem:rec}
Given a tree $\Gamma$, a cutset $C$, interactions $\J$, boundary condition
$\delta$ and $v\in C^i$, 
let $\{ w_1 , \ldots , w_k \}$ be the children of $v$.  Let $J_1 , \ldots , 
J_k$ denote the values of $\J (v , w_1) , \ldots , \J (v , w_k)$.  Then 
\begin{equation} \label{eq:recurse}
f_{C,v}^{\J,\dd} = \poi (\K_{J_1} f^{\J , \delta}_{C , w_1} , 
   \ldots , \K_{J_k} f^{\J , \delta}_{C , w_k}) \, ,
\end{equation} 
where when $w_i\in C$, $f^{\J , \delta}_{C , w_i}$ is taken to be the point
mass at $\delta(w_i)$ and convention~(\ref{eq:pointconv}) is in effect.
\end{lem}

\noindent{\bf Proof.}
Passing to the subtree $\Gamma (v)$, we may assume without loss of
generality that $v = o$.  Also assume without loss of generality that
$w_1 , \ldots , w_k$ are numbered so that
for some $s$, $w_i \in C^i$ for $i \leq s$ and $w_i \in C$ for $i > s$.
For $i\le s$, let $C(w_i) = C \cap \Gamma (w_i)$.  For such $i$, by 
definition, $f_i := f^{\J , \dd}_{C,w_i}$ is the marginal at $w_i$ 
of the measure $\mu_i := \mu^{\J, \dd}_{C(w_i) , w_i}$ on 
configurations on $\Gamma (w_i)\cap C^i$, where $\J$ and $\p$ 
are restricted to $E (\Gamma (w_i))$ and $C(w_i)$ respectively.
Let $\Gamma_r$ denote the induced subgraph of $\Gamma$ whose vertices 
are the union of $\{ o \}$, $\Gamma (w_1) , \ldots , \Gamma (w_r)$.
We prove by induction on $r$ that the density $g_r$ at the root of
$\Gamma_r$ of the analogue of $\mu^{\J , \delta}_C$ for $\Gamma_r$ is equal to
$$\poi (\K_{J_1} f^{\J , \delta}_{C , w_1} , \ldots , 
   \K_{J_r} f^{\J , \delta}_{C , w_r}) \, ;$$
The case $r = k$ is the desired conclusion.

To prove the $r=1$ step, use~(\ref{eq:fdecomp}) with $v_1 = o$, 
$A_1 = \{ o \}$, $C_1 = \emptyset$, $v_2 = w_1$, $A_2 = \Gamma (w_1)$ 
and $C_2 = C(w_1)$.  If $w_1 \in C$, the $r=1$ case is 
trivially true, so assume $s \geq 1$.
The measure $\mu_1$ is uniform on $\S$ since $C(v) = \emptyset$.
Thus from~(\ref{eq:fdecomp}) we find that
$$g_1 (y) = c \int e^{-J_1 H (y,z)} f_1 (z) \, dz 
   = (\K_{J_1} f_1) (y)$$
which proves the $r = 1$ case.  

For $1 < r \leq s$, use~(\ref{eq:fdecomp}) with $A_1 = \Gamma_{r-1}$,
$v_1 = o$, $C_1 = \Gamma_{r-1} \cap C$, $A_2 = \Gamma (w_r)$,
$v_2 = w_r$ and $C_2 = \Gamma (w_r) \cap C$.  
Using~(\ref{eq:fdecomp}) we find that
\begin{eqnarray*}
g_r (y) & = & c g_{r-1} (y) \int e^{-J_r H (y,z)} 
   f_r (z) \, \dx (z) \\[1ex]
& = & c g_{r-1} (y) (\K_{J_r} f_r) (y) \\[1ex]
& = & (\poi (g_{r-1} , \K_{J_r} f_r)) (y) \, .
\end{eqnarray*}
By associativity of $\poi$ the induction step is completed for $r \leq s$.

Finally, if $r > s$, then the difference between $H (\eta)$ on
$\Gamma_{r-1}$ and $H (\eta)$ on $\Gamma_r$ is just $- J_r
H(\eta (o) , \delta (w_r))$, so 
$$g_r (y) = c g_{r-1} (y) \exp (- J_r 
H(y , \delta (w_r))) = \left ( \poi 
   (g_{r-1} , \K_{J_r} f_r) \right ) (y)$$
by the convention~(\ref{eq:pointconv}), and associativity of $\poi$
completes the induction as before.
$\qed$

Another consequence of Lemma~\ref{lem:decomps} is Lemma~\ref{lem:free},
giving the existence of a natural and well defined free boundary measure.

\noindent{\bf Proof of Lemma~\protect{\ref{lem:free}}.}
Observe that in~(\ref{eq:mudecomp2}), if $f_2 \equiv 1$ then 
the integral against $z$ is independent of $y$, so one has 
$\mu^{(1)} = \mu_1$.  Let $F$ be any cutset and $w \in F^i$ be
chosen so each of its children $v_1 , \ldots , v_k$ is in $F$.  
Applying our observation inductively to eliminate each child 
of $w$ in turn, we see that the projection of $\mufree_F$ onto 
$\S^{F^i \setminus \{ w \}}$ is just $\mufree_{F'}$ where 
$F' = F \cup \{ w \} \setminus \{ v_1 , \ldots , v_k \}$.  

Given cutsets $C$ and $D$ with $D \cap C^i \neq \emptyset$, choose 
$v \in D \cap C^i$ and $w \geq v$ maximal in $C^i$.  Then all
children of $w$ are in $C$.  Applying the previous paragraph
with $F = C$, we see that $\mufree_C$ agrees with $\mufree_{F'}$.
Continually reducing in this way, we conclude that on $C^i \cap D^i$
$\mufree_C$ agrees with $\mufree_Q$ where $Q$ is the
exterior boundary of $C^i \cap D^i$.  The same argument shows
that $\mufree_D$ agrees with $\mufree_Q$, which finishes the
proof of the lemma.
$\qed$

According to Lemma~\ref{lem:rec}, if, for $J>0$, we define $\PP (J)$ to be 
the smallest class of densities containing each $K_{J' , y}$ for 
$J' \in (0,J]$ and $y \in \S$ and closed under $\K_{J'}$ for $J' \in (0,J]$
and $\poi$, then, when $\J$ is strictly positive and
bounded by $J$, each density $f^{\J , \p}_{C , v}$ is an element of $\PP(J)$. 
Similarly, if $\PP_+(J)$ is taken to be 
the smallest class of densities containing each $K_{J'}$ for 
$J' \in (0,J]$ and closed under $\K_{J'}$ for $J' \in (0,J]$
and $\poi$, then, when $\J$ is strictly positive and 
bounded by $J$, each density 
$f^{\J , +}_{C , v}$ is an element of $\PP_+(J)$. We also let
$\PP:=\bigcup_{J> 0}\PP(J)$ and $\PP_+:=\bigcup_{J> 0}\PP_+(J)$.

This leads to the following lemma whose proof is left to the reader.   

\begin{lem} \label{lem:unifbd} 
Suppose the interaction strengths $\{ \J (e) \}$ are bounded above
by some constant.  
Then there exist constants $0 < B_{\rm min} < B_{\rm max}$ such that 
for every $C , \delta$ and $v \in C^i$, the one-dimensional marginal
of $\mu^{\delta}_C$ at $v$ is absolutely
continuous with respect to $\dx$ with a
density function in $[B_{\rm min} , B_{\rm max}]$.
It follows, since the above properties are closed under convex combinations,
that all one-dimensional marginals of any Gibbs state have densities
in $[B_{\rm min} , B_{\rm max}]$. Similarly, the $k$-dimensional marginals have
densities in the interval $[B_{\rm min}^{(k)} , B_{\rm max}^{(k)}]$
for some constants $0< B_{\rm min}^{(k)} < B_{\rm max}^{(k)}$.
In addition, the family of all one--dimensional densities which arise as
above is an equicontinuous family.
\end{lem}

The usefulness of the equicontinuity property is that the following
easily proved lemma (whose proof is also left to the reader) tells us that
in determining weak convergence to $\dx$, it is equivalent to look to see if
there is convergence in $L^\infty$ of the associated densities to 1.
\begin{lem} \label{lem:converge}
Let $(X,d)$ be a compact metric space and $\mu$ a probability measure on
$X$ with full support. If $\{f_n\}$ is an equicontinuous family of
probability densities (with respect to $\mu$), then
$$
\lim_{n\to\infty} \|f_n-1\|_{\infty} = 0 \mbox{ if and only if } 
\lim_{n\to\infty} f_n d\mu = \mu \mbox{ weakly }.
$$
\end{lem}

Using this, we can prove the equivalence of phase transition and
symmetry breaking on trees (Proposition~\ref{prop:SB=}).

\noindent{\bf Proof of Proposition~\ref{prop:SB=}.}
(i) implies (ii) is trivial. For (ii) implying (iii), assume we have
a vertex $v$, a sequence of cutsets $C_n\to\infty$ and boundary 
conditions $\p_n$ on $C_n$ such that
$$
\inf_n \|f_{C_n,v}^{\p_n}-1\|_\infty \neq 0.
$$
Clearly we obtain the same result if we change $\p_n$ on 
$C_n\setminus \Gamma(v)$ to anything, in particular, if we take no (i.e., 
free) boundary condition there. We then take any weak limit of these 
measures as $n\to\infty$. This will yield a Gibbs state and by the first line
of the proof of Lemma~\ref{lem:free}, together with 
Lemma~\ref{lem:converge},
the marginal density at $v$ of this Gibbs state is not 1, which proves (iii).
(iii) implies (iv) is also trivial of course. To see that (iv) implies (i),
note that if there is PT, then there exists an extremal Gibbs state 
$\mu\neq \mufree$. Choose a cutset $C$ such that
$\mu\neq \mufree$ when restricted to $C^i$. If (i) fails, then for all
$v\in C$, there exists a
sequence of cutsets $C_n\to\infty$ such that for all boundary 
conditions $\p_n$ on $C_n$ we have that
\begin{equation} \label{eq:ivgivesi}
\inf_n \|f_{C_n,v}^{\p_n}-1\|_\infty = 0.
\end{equation} 
Clearly, because of the geometry, $\{C_n\}$ can be chosen independent of $v$. 
Since $\mu$ is extremal, it is known (see Theorem 7.12(b) in \cite{Ge}, p. 122)
that there exist boundary conditions
$\p_n'$ on $C_n$ so that $\mu_{C_n}^{\p_n'} \rightarrow \mu$ weakly.
However, by (\ref{eq:ivgivesi}) and Lemma~\ref{lem:rec}, $\mu$ must equal
$\mufree$ on $C^i$, a contradiction. $\qed$

\subsection{Distance regular spaces} \label{sub:drs}

Our primary interest in this paper is in the Heisenberg models.
Nevertheless, it turns out that many of the properties of the Heisenberg model
hold in the more general context of distance regular spaces.
A {\bf distance regular graph} is a finite graph for
which the size of the set $\{ z : d(x,z) = a , d (y,z) = b \}$ 
depends on $x$ and $y$ only through the value of $d(x,y)$ where
$d(x,y)$ is the usual graph distance between $x$ and $y$.
We generalize this by saying that the metric space $(\S,d)$ 
with probability measure 
$\dx$ is {\bf distance regular} if the law of the pair $(d(x,Z) , d(y,Z))$
when $Z$ has law $\dx$ depends only on $d(x,y)$.  
In particular, when the action of $G$ on $\S$ is distance transitive
(in addition to preserving $d$ and $\dx$), 
meaning that $(x,y)$ can be mapped to any $(x' , y')$ with 
$d(x,y) = d(x' , y')$, it follows easily that $(\S,d, \dx)$ is 
distance regular.  All the examples we have mentioned so far are distance 
transitive (and hence distance regular) 
except for the rotor model which is still distance regular.
(For an example of a graph showing that the full automorphism group
acting distance transitively
is strictly stronger than the assumption of distance 
regularity, see~\cite{AVLF} or 
{\it Additional Result} {\bf 23b} of~\cite{Big}.)  

We present some of the background in this generality not because 
we are fond of gratuitous generalization but because we find the 
reasoning clearer, and because it seems reasonable that
someone in the future might study a particle system whose 
spin states are elements of some distance regular space, 
such as real projective space or the discrete $n$-cube.  The primary
consequence of distance regularity is that it allows one to define
a commutative convolution on a certain subspace of $L^2$.

\begin{defn} 
Let $\lts$ denote the space $L^2 (\dx)$, and let $\ltsg$ denote
the space of functions $f \in \lts$ for which $f(x)$ depends
only on $d(x , \zero)$.  For $f \in \ltsg$, define a function $\fb$
on $\{d(\zero,y)\}_{y\in \S}$ 
by $\fb (r) := f(x)$ where $x$ is such that $d(\zero,x) = r$.
\end{defn}

\begin{defn}
If $(\S,\dx)$ is distance regular,
define a commutative convolution operation on $\ltsg \times \ltsg$ by 
$$f * h (x) :=  \int_{\S} h(y) \fb (d(x,y)) \, \dx (y) = 
   \int_{[0,\infty)^2} \fb (u) \hb (v) \, d\pi_x(u,v)$$
where $\pi_x$ is the law of $(d(x,Z) , d(\zero , Z))$ for a variable
$Z$ with law $\dx$.  It is clear from the definition of a
distance regular space that $(d(x,Z) , d(\zero,Z))$ and
$(d(\zero,Z) , d(x,Z))$ are equal in distribution implying
that $f * h =h * f$ and that, since $\pi_x$ only depends on $d(x,\zero)$,
$f,h\in \ltsg$ implies that $f * h \in \ltsg$.
\end{defn} 

The following lemma is straightforward and left to the reader.

\begin{lem} \label{lem:dt}
For all $J\ge 0$, $K_J\in \ltsg$ and for all $h\in \lts$,
$\K_J(h)(x)$ (defined in~(\ref{eq:conv})) is equal to 
$\int_{\S} h(y) \overline{K_J} (d(x,y)) \, \dx (y)$.
In particular, if $(\S , \dx)$ is distance regular, 
then the operators $\K_J$ map $\ltsg$ into itself and 
$\K_J(h) =K_J * h$ for all $h \in \ltsg$.
\end{lem}

We believe that for most distance regular spaces, one can verify the 
necessary hypotheses of Theorems~\ref{th:gen ii} and~\ref{th:gen i} 
below in the same way as
we will do for the Heisenberg models in detail 
in the next section. Doing this however
would take us too far afield and so we content ourselves with pointing out
to the reader that much of this probably can be done, and after analyzing
the Heisenberg models in Section~\ref{sub:sphere}, 
explain how to carry much of this out in the context of distance regular 
graphs in Section~\ref{sub:finite}.  

\subsection{Heisenberg models} \label{sub:sphere}

In this subsection, we consider Example \ref{eg:spherical} in  
Section~\ref{sec:one}
and so we have $\S = S^d$, $d \geq 1$, the unit sphere in
$(d+1)$--dimensional Euclidean space with the corresponding
$G, d, \dx$ and $H$.
Recall that this is distance transitive for $d\ge 2$ (and hence 
distance regular) and distance regular for $d=1$.
The following lemma allows us to set up coordinates in which our bookkeeping
will be manageable. It is certainly well known.

\begin{lem} \label{lem:spherical} For any $d\ge 1$,
there exist real--valued functions $\f_0 , \f_1 , \f_2 , \ldots \in \ltsg$
$(\S=S^d)$,
orthogonal under the inner product 
$\langle f,g \rangle = \int_{\S} f \overline{g} \, \dx$, 
such that $\f_n$ is a polynomial of degree exactly $n$ in $x\cdot \zero$,
and such that the following properties hold. \\
(1) $\f_0 (x) \equiv 1$ and $\f_1 (x) = x \cdot \zero$. \\
(2) $1 = \f_j (\zero) = \sup_{x \in \S} |\f_j (x)|$, for all $j$. \\
(3) $\f_i \f_j = \sum_{r\ge 0} q^r_{ij} \f_r$, where the coefficients 
   $q^r_{ij}$ are nonnegative and $\sum_r q^r_{ij} = 1$. \\
(4) $\f_i * \f_j = \gamma_j \delta_{ij} \f_j$, where 
   $\gamma_j := \f_j * \f_j (\zero) = \int \f_j^2(x) \, \dx(x)$.  \\
(5) The functions $\f_j$ are eigenfunctions of any convolution
   operator, that is, $f * \f_j = c \f_j$ for any $f \in \ltsg$. \\
(6) Any $f \in \ltsg$ can be written as a convergent series 
   $f(x) = \sum_{j\ge 0} a_j (f) \f_j(x)$ (in the $L^2$ sense), 
where the complex numbers $a_j(f)$ 
are given by $a_j(f) : = \gamma_j^{-1} \int f(x) \f_j (x) \, \dx (x) .$ \\
(7) For $f,g \in \ltsg$, we have $a_j(f * g)= \gamma_j a_j(f) a_j(g)$.
%%\item The eigenvalues of the convolution operator $g \mapsto f * g$
%%   are precisely $\gamma_j a_j(f)$ which is equal to $\int f \f_j \, \dx$.
\end{lem}

\noindent{\bf Proof.}  For each $\alpha, \beta >-1$, define the 
Jacobi polynomials $\{\P^{(\alpha , \beta)}_n(r)\}_{n\ge 0}$ by
\begin{equation} \label{eq:rod}
(1 - r)^\alpha (1 + r)^\beta \P_n^{(\alpha , \beta)} (r) = 
   {(-1)^n \over 2^n n!} {d^n \over dr^n} \left [ (1 - r)^{n + \alpha}
   (1 + r)^{n + \beta} \right ] \, .
\end{equation}
(The Jacobi polynomials are usually defined differently in which
case~(\ref{eq:rod}) becomes what is known as Rodrigues' formula but we shall 
use~(\ref{eq:rod}) as our definition; when $\alpha=\beta$, which is the case 
relevant to us, these are the ultraspherical polynomials.)

For any given $d\ge 1$, we let, for $n\ge 0$,
$$
\f_n(x):=
{\P^{({d \over 2}-1 ,{d \over 2}-1)}_n (x\cdot \zero)
\over
\P^{({d \over 2}-1 ,{d \over 2}-1)}_n (1)}.
$$
By p.254 in \cite{R}, $\P^{(\alpha , \beta)}_n$ is a polynomial  of degree
exactly $n$.
By p.259 in \cite{R}, the collection $\{\P^{(\alpha , \beta)}_n\}_{n\ge 0}$ 
are orthogonal on $[-1,1]$ with respect to the weight function
$(1-r)^\alpha (1+r)^\beta$. A change of variables then shows that the $\f_n$'s
are orthogonal in $\lts$. 

(1) is then an easy calculation, the first equality in (2) is trivial 
while the second equality is in \cite{R}, p.278 and 281. 
(3) is in~\cite{Askey74}, 
p.41. (4) and (5) follow from the Funk--Hecke Theorem (\cite{N}, p.195)
(the calculation of $\gamma_j$ being trivial).
Since the subspace generated by the 
$\{\P^{({d \over 2}-1 ,{d \over 2}-1)}_n(r)\}$'s are uniformly dense 
in $C([-1,1])$ by the Stone-Weierstrass Theorem, it easily follows that the 
subspace generated by the $\f_n$'s are uniformly dense in $\ltsg\cap C(\S)$.
Hence the $\f_n$'s are a basis for $\ltsg$ and (6) follows.
Finally, (4) and (6) together yield (7). $\qed$

Note that for all $f,g\in \ltsg$, we have that
$fg\in \ltsg$ provided $fg\in \lts$.
Since $\f_n$ is a polynomial of degree exactly $n$ in $x\cdot \zero$,
the greatest $r$ for which
$q^r_{ij} \neq 0$ must be $i + j$.  From this and the nonnegativity 
of the $q^r_{ij}$'s, it follows that for $\lambda > 0$ the function 
$e^{\lambda \f_1(x)} = \sum_{n \geq 0} \lambda^n \f_1 (x)^n / n!$ has
\begin{equation} \label{eq:viii}
a_j (e^{\lambda \f_1}) > 0 , \mbox{  for all } j\ge 0.
\end{equation}
It follows from Lemmas~\ref{lem:rec}, \ref{lem:dt} 
and \ref{lem:spherical}(3,4) that $\PP_+\subseteq \ltsg$ and that
for all $g \in \PP_+$,
\begin{equation} \label{eq:ix}
a_j (g) > 0, \mbox{  for all } j\ge 0.
\end{equation}

\begin{defn}
Define the $A$ norm on $\ltsg$ by
$$||f||_A = \sum_{j\ge 0} |a_j (f)| ,$$
provided it is finite.
\end{defn}

 From the fact that $\sum_{r\ge 0} q^r_{ij} = 1$, one can easily show 
that for all $f,g\in\ltsg$ with $fg\in\ltsg$,
\begin{equation} \label{eq:submult}
||fg||_A \leq ||f||_A ||g||_A , 
\end{equation}
and that equality holds if $f , g \in \PP_+$. An easy computation also 
shows that $||e^{\lambda \f_1(x)}||_A =e^{\lambda} <\infty$ 
for all $\lambda \ge 0$
and hence by Lemmas~\ref{lem:rec} and~\ref{lem:spherical}(4)
and~(\ref{eq:submult}), $||f||_A<\infty$ for all $f\in\PP_+$.
Also, it follows from~(\ref{eq:ix}), Lemma~\ref{lem:spherical}(2,6),
the fact that $\int f \, \dx = 1$ for all $f \in \PP_+$ and the fact
that $\PP_+\subseteq \ltsg$ that for $f \in \PP_+$,
\begin{equation} \label{eq:x}
1 + ||f - 1||_A = ||f||_A = f (\zero) = ||f||_\infty = 
   1 + ||f - 1||_\infty . 
\end{equation}
The last equality is obtained by observing that $\le$ is clear while
$ ||g||_\infty \le ||g ||_A $ for all $g\in \ltsg$ is also clear.

\begin{lem} \label{lem:taylor}
There exists a function $o$ with $\displaystyle{\lim_{h \to 0} {o(h) \over h} 
= 0}$ such that for all $h_1 , \ldots , h_k \in \PP_+$
with $k \leq B$,
\begin{equation} \label{eq:xi}
|| \poi (h_1 , \ldots , h_k) - 1 - \sum_{i=1}^k (h_i - 1)||_A
   \leq o(\max_i ||h_i - 1||_A) ,
\end{equation}
provided $\max_i ||h_i - 1||_A\le 1$.
\end{lem}

\noindent{\bf Proof.} Write
\begin{equation} \label{eq:new01}
|| \prod_{i=1}^k h_i - 1 - \sum_{i=1}^k (h_i - 1)||_A = 
|| \sum_{{A\subseteq\{1,\ldots,k\}\atop |A|\ge 2}} \prod_{i\in A}(h_i-1)||_A.
\end{equation}
Then $\max_i ||h_i - 1|| \leq 1$ and
submultiplicativity~(\ref{eq:submult})
of $|| \cdot ||_A$ implies this is at most
$$ 2^k (\max_i ||h_i - 1||_A)^2 .$$
Next, since $\int (h_i - 1) \, \dx = 0$ for $1 \leq i \leq k$, 
we similarly obtain
$$\left|\int \prod_{i=1}^k h_i - 1\right| \leq 2^k (\max_i ||h_i - 1||_A)^2.$$
We then have
$$
|| \poi (h_1 , \ldots , h_k) - \prod_{i=1}^k h_i||_A =
{1\over \int \prod_{i=1}^k h_i}\left|\int \prod_{i=1}^k h_i - 1\right|
|| \prod_{i=1}^k h_i ||_A
\leq 4^k (\max_i ||h_i - 1||_A)^2 ,$$
since $|| \prod_{i=1}^k h_i ||_A\le 2^k$ and
$\int \prod_{i=1}^k h_i\ge 1$ by 
the positivity of the $q^r_{ij}$ and~(\ref{eq:ix}).
A use of the triangle inequality completes the proof. $\qed$

We note five facts that follow easily from the above, but
which will be useful later on in generalizing our results.
%%Theorems~\ref{th:gen ii} and~\ref{th:gen i}. 
Let $\ssp$ be the linear subspace of $\ltsg$ spanned by $\PP_+$,
$\sspj$ be the linear subspace of $\ltsg$ spanned by $\PP_+(J)$
and $|| \calk_{J'} ||_A$ denote the operator norm of $\calk_{J'}$ on
$(\ssp,||\,\,||_A)$.
\begin{equation} \label{eq:xii}
\lim_{J' \to 0} ||K_{J'} - 1||_A = 0 ;
\end{equation}

\begin{equation} \label{eq:present iii}
c_1:=\sup_{f\in\ssp ,f \neq 1} {||f - 1||_\infty \over ||f - 1||_A} < \infty ;
\end{equation}

\begin{equation} \label{eq:present iii'}
c_2:=\inf_{f\in\PP_+,f\neq 1} {||f - 1||_\infty \over ||f - 1||_A} > 0 ;
\end{equation}

\begin{equation} \label{eq:xiii}
\mbox{ For all } J' \geq 0, \,\,\, || \calk_{J'} ||_A \leq 1;
\end{equation}

There exist $a,b\in \S$ such that for all  $f\in \PP_+$,
\begin{equation} \label{eq:xiiii}
f(a)=\sup_{x\in\S} f(x) \mbox{ and } f(b)=\inf_{x\in\S} f(x).
\end{equation}

(\ref{eq:xiii}), for example, follows immediately from Lemmas~\ref{lem:dt} 
and~\ref{lem:spherical}(7) and the fact that $|\gamma_n a_n(g)|\le 1$ for any 
probability density function $g\in \ltsg$.

The results on Heisenberg models presented thus far are parallel to 
the results obtainable for any finite distance regular graph (see
the next subsection).  One useful result that is not true for general 
distance regular models depends on the following obvious geometric property 
of the sphere:
$$|\{ z : d(x,z) \leq a , d(y,z) \leq b \}|$$
is a nonincreasing function of $d(x,y)$ for any fixed $a$ and $b$ where
$|\,\,|$ denotes surface measure.
[Proof: For $S^1$, this is obvious. For $S^d$, $d\ge 2$, by symmetry, we can
assume that $x=(0,\ldots,0,1)$ and $y=(\cos\theta,0,\ldots, 0, \sin\theta)$
(both vectors with $d+1$ coordinates). Write $S^d$ as
$$
\cup_{u\in [-1,1]^{d-1}} A_u
$$
where 
$$
A_u:=S^d\cap\{(a_1,\ldots,a_{d+1}):(a_2,\ldots,a_{d})=u\}.
$$
Each $A_u$ is a circle (or is empty) and so essentially by the 1--dimensional
case, we have the desired behaviour on each $A_u$ (using 1--dimensional
Lebesgue measure) and by Fubini's Theorem, we obtain the desired result
on $S^d$.]

Calling a function $f\in\ltsg$ nonincreasing if the corresponding $\fb$ is
nonincreasing, the latter can be seen to be 
equivalent to the property that $\one_{d(x , \zero) \leq a} * 
\one_{d(x , \zero) \leq b}$
is nonincreasing, and by taking linear combinations, this is equivalent 
to $f * g$ being nonincreasing for all nonincreasing $f$ and $g$ in $\ltsg$.  
Since $K_J$ is nonincreasing for all $J$,
it follows from the fundamental recursion that
\begin{equation} \label{eq:xiv}
f \in \PP_+ \Rightarrow f \mbox{ is nonincreasing}  .
\end{equation}

\begin{lem} \label{lem:incr}
For any positive nonincreasing $f \in \ltsg$, 
$$
\left|\int_{\S} f \f_n \, \dx\right| \le \int_{\S} f \f_1 \, \dx
$$
for all $n \ge 1 $.
\end{lem}

\noindent{\bf Proof.}  It suffices to prove this for functions of the form
$f(x) = \one_{\{x \cdot \zero \geq t\}}$ with $t\in [-1,1]$.  
We rely on explicit formulae for the functions $\{ \f_n \}$. 
Letting $\alpha = d/2 - 1$, a change of variables yields
$$\int_{\S} f \f_n \, \dx 
= s_d^{-1}\int_t^1 {P_n^{(\alpha , \alpha)} (r) \over P_n^{(\alpha ,
   \alpha)} (1)} (1 - r^2)^\alpha \, dr,$$ 
where 
$$
s_d=\int_{-1}^1 (1-r^2)^\alpha dr
$$
and
$\P_n^{(\alpha , \alpha)}$ is the Jacobi polynomial defined earlier.

Taking the indefinite integral of each side in~(\ref{eq:rod})
with $\beta = \alpha$ yields
\begin{eqnarray*}
\int (1 - r)^\alpha (1 + r)^\alpha P_n^{(\alpha , \alpha)} (r) \, dr & = &
   {(-1)^n \over 2^n n!} {d^{n-1} \over dr^{n-1}} \left [ (1 - r)^{n + \alpha}
   (1 + r)^{n + \alpha} \right ] \\[2ex]
& = & {- 1 \over 2n} (1 - r^2)^{\alpha + 1} P_{n-1}^{(\alpha + 1 , 
   \alpha + 1)} (r)\, .
\end{eqnarray*}
Evaluating at 1 and $t$ gives
\begin{eqnarray*}
\int_{\S} f \f_n \, \dx 
 & = & s_d^{-1}\int_t^1 {P_n^{(\alpha , \alpha)} (r) (1 - r^2)^\alpha \over
   P_n^{(\alpha , \alpha)} (1)} \, dr \\[2ex]
& = & s_d^{-1}
{P_{n-1}^{(\alpha + 1 , \alpha + 1)} (t) (1 - t^2)^{\alpha + 1} \over
   2 n P_n^{(\alpha , \alpha)} (1)} .
\end{eqnarray*}
When $n = 1$, using~(\ref{eq:rod}),
this is just $s_d^{-1}(1 - t^2)^{\alpha + 1} / 2(1+\alpha)$.  Dividing, we get 
$${ \int_{\S} f \f_n \, \dx \over \int_{\S} f \f_1 \, \dx } 
= {P_{n-1}^{(\alpha + 1 , \alpha + 1)} (t) (1+\alpha)
   \over n P_n^{(\alpha , \alpha)} (1)} = 
   {P_{n-1}^{(\alpha + 1 , \alpha + 1)} (t) \over
   P_{n-1}^{(\alpha + 1 , \alpha + 1)} (1)} \cdot 
   {P_{n-1}^{(\alpha + 1 , \alpha + 1)} (1) \over
   n P_n^{(\alpha , \alpha)} (1)}\cdot(1+\alpha).$$
The first term in the product is bounded in absolute value by 1.  
By \cite{Askey74}, p.7,
$$P_n^{(\alpha , \alpha)} (1) = {\alpha + n \choose n} ,$$
and so we see that the second term is $1 / (\alpha + 1)$, completing the
proof of the lemma. $\qed$

\noindent{\em Remark:}
The case $d = 1$ can also be handled by a rearrangement lemma.

\begin{defn} \label{defn:op}
Define a linear functional $L$ on $\ltsg$ by 
$L (g):= \int_{\S}g(x)\f_1(x) \dx(x)$ $(= \gamma_1 a_1 (g))$
and set $\op = L(K_J)$. (Recall that $\f_1,\gamma_1$ and $a_1$ are defined
in Lemma~\ref{lem:spherical}.)
\end{defn}

It follows from Lemmas~\ref{lem:dt},~\ref{lem:spherical}(7) 
and~\ref{lem:incr},~(\ref{eq:xiv}) and an easy computation  that
\begin{equation} \label{eq:ixix}
||\calk_J f-1||_A\le \op ||f-1||_A \mbox{  for all  } f\in\PP_+(J).
\end{equation}
%%that $\op$ is the operator norm of $\K_J$ on 
%%$\{f\in\ltsg:\, ||f||_A <\infty, \int_{\S} f \dx=0\}$.
In the following inequalities, we denote 
$\rho := \op$.  For $f \in \PP_+(J)$, it also follows easily that
\begin{equation} \label{eq:(a)}
L (\K_J f - 1 ) \geq \rho L (f - 1)
\end{equation}
and that there is a constant $c_3$ such that for all $f\in\,\,\sspj$,
\begin{equation} \label{eq:(c)}
|L(f)| \leq c_3 ||f||_A .
\end{equation}
(We can of course take $c_3$ to be 1, but we leave the condition written 
in this more general form for use as a hypothesis in Theorem~\ref{th:gen i}.)

Putting together the results of Lemmas~\ref{lem:spherical} and~\ref{lem:incr},
as well as~(\ref{eq:viii}),~(\ref{eq:ix}) and~(\ref{eq:xiv}),
gives the following corollary.

\begin{cor} \label{cor:(b)}
For all $J\ge 0$, 
there is a constant $c_4 > 0$ such that for all $f \in \PP_+(J)$,
\begin{equation} \label{eq:(b)}
L (f) \geq c_4 ||f - 1||_A .
\end{equation}
\end{cor}

\noindent{\bf Proof.}  Fix $f \in \PP_+(J)$.  If $f = K_{J'}$ for
some $J' \in (0, J]$, we argue as follows. As 
$||e^{\lambda \f_1(x)}||_A =e^{\lambda}$ (which we mentioned earlier)
and $K_{J'}(x)=e^{J'\f_1(x)}/\int e^{J'\f_1(x)}\dx(x)$, we have
\begin{eqnarray*}
||K_{J'} - 1||_A  & = & ||K_{J'} ||_A -1 
\\[2ex]
& = & {e^{J'}\over \int e^{J'\f_1(x)}\dx(x)} -1
\\[2ex]
& \le & e^{2J'} -1.
\end{eqnarray*}
Next, 
\begin{eqnarray*}
L (K_{J'}) & = & {1\over \int e^{J'\f_1(x)}\dx(x)} 
\int e^{J'\f_1(x)}\f_1(x)\dx(x) \\[2ex]
& = & {1\over \int e^{J'\f_1(x)}\dx(x)} 
\sum_{k=0}^\infty {(J')^k\over k!} \int \f_1^{k+1}(x)\dx(x).
\end{eqnarray*}
By Lemma~\ref{lem:spherical}(3), all terms in the sum are nonnegative and
by Lemma~\ref{lem:spherical}(4), the $k=1$ term is $J'\gamma_1$.
Hence $L (K_{J'})\ge J'\gamma_1/e^{J'}$. 
Since
$$
\inf_{J'\in (0,J]} {J'\gamma_1 \over  e^{J'}(e^{2J'} -1)  } >0,
$$
we can find a $c_4$ in this case.
%%the result follows from~(\ref{eq:viii}) and compactness; by direct 
%%computation we can take $\eps \to 0$ and retain a uniform lower 
%%bound on $L(f)/||f - 1||_A$.  

Otherwise, by the fundamental recursion, we may represent $f$ as
$\poi (\K_{J_1} h_1 , \ldots , \K_{J_k} h_k)$ with each
$h_i$ either in $\PP_+(J)$ or equal to $\delta_{\zero}$ and each 
$J_i\in (0,J]$. Define $g_i = \K_{J_i} h_i - 1$.
Let $m :=\inf_{0< J'\le J} a_1 (K_{J'}) / \sum_{n > 0} a_n (K_{J'})$ which is
strictly positive by the above. It follows that if %%$J_i >0$ and 
$h_i\in \PP_+(J)$ (the case $h_i=\delta_{\zero}$ is already done),
$${L (g_i) \over ||g_i ||_A} = {a_1 (K_{J_i}) a_1 (h_i) \gamma_1^2\over
   \sum_{n > 0} a_n (K_{J_i}) a_n (h_i) \gamma_n} \geq m\gamma_1$$
by Lemma~\ref{lem:spherical}(7) and
since $a_1(h_i)\gamma_1 \ge a_n(h_i)\gamma_n$ for all $n\ge 1$ by
Lemma~\ref{lem:incr} and~(\ref{eq:xiv}).
Let $h = \prod_{i=1}^k \K_{J_i} h_i$.  Then
$L (h) = L(1 + \sum_{i=1}^k g_i + Q)$, where $Q$ is a sum 
of monomials in $\{ g_i \}$. Using $q^r_{ij} \ge 0$
and~(\ref{eq:ix}), we have that $L(Q) \geq 0$, and hence
\begin{equation} \label{eq:new02}
L(h) \geq \sum_{i=1}^k L(g_i) \geq m \gamma_1\sum_{i=1}^k ||g_i ||_A .
\end{equation}
On the other hand, for any $B$ and $M$, there is $C = C(M,B)$ such that
if $x_1 , \ldots , x_k \in (0,M)$ with $k\le B$, then
$$ -1 + \prod_{i=1}^k (1 + x_i) \leq C \sum_{i=1}^k x_i .$$
Next, 
the positivity of the $q^r_{ij}$ implies 
$\int_{\S} h(x) \, \dx(x) = a_0 (h) \geq 1$. It follows that 
$$||h - 1||_A =  -1+ ||h ||_A = 
-1 + \prod_{i=1}^k ||g_i + 1||_A \leq C \sum_{i=1}^k 
   ||g_i||_A $$
for some constant $C$ since $||g_i+1 ||_A=||g_i ||_A +1$ and
$||g_i+1 ||_A$ clearly has a universal upper bound. [To see the latter
statement, one notes that
$$
\sup_{0<J' \le J} ||K_{J'}||_A <\infty,
$$
$$
||K_{J'} *f||_A \le ||K_{J'} ||_A
$$
for any probability density function $f\in \ltsg$ 
(by Lemma~\ref{lem:spherical}(7)),~(\ref{eq:submult}) and the fact that we 
never have more than $B$ terms in our pointwise products imply that
$$
\sup_{f\in \PP_+(J)}||f||_A \le 
\left(\sup_{0<J' \le J} ||K_{J'}||_A \right)^B <\infty.]
$$
Putting this together with~(\ref{eq:new02}) gives
$${L(h) \over ||h - 1||_A} \geq {m \gamma_1\over C} \, .$$

Finally, letting $f = h / \left(\int_{\S} h(x) \, \dx(x)\right)$, we obtain
\begin{eqnarray*}
||h - 1||_A  & \ge &  \sum_{n\ge 1} a_n(h) 
\\[2ex]
& = & \sum_{n\ge 1} \left[\int_{\S} h(x)\dx(x)\right] a_n(f) \\[2ex]
& = &  \left[\int_{\S} h(x)\dx(x)\right]  ||f - 1||_A.
\end{eqnarray*}
Hence
$$
{L(f)\over ||f - 1||_A} \ge {L(h)\over ||h - 1||_A} \ge {m \gamma_1\over C}
$$
and we're done. $\qed$

\subsection{Distance regular graphs} \label{sub:finite}

For the remainder of this section, we suppose that $\S$ is the vertex
set of a finite, connected, distance regular graph, that $d(x,y)$ is 
the graph distance, and that the energy $H (x,y)$ depends only 
on $d(x,y)$.  The Potts models fit into this 
framework, with the respective graphs being 
the complete graph $K_q$ on $q$ vertices.  
All the results we need follow
in fact from an even weaker assumption, namely that $\S$ is an
{\it association scheme}.  For the definition of association schemes
and the proofs of the relevant results, see~\cite{BCN} or~\cite{Ter98}.
By developing 
the analogue of Lemma~\ref{lem:spherical} for distance regular graphs, 
we will illustrate the extent to which our results are independent of
the special properties of the Heisenberg model.

We have a distinguished element $\zero \in \S$ and the measure $\dx$ will of
course be normalized counting measure $|\S|^{-1} \sum_{x \in \S} \delta_x$.
The spaces $\lts$ and $\ltsg$ are then simply finite dimensional vector
spaces with respective dimensions $|\S|$ and $1+D$, where $D$ is the 
diameter of the graph $\S$.  

Denote by $M(\S)$ the space of matrices with rows and columns
indexed by $\S$, thought of as linear maps from $\lts$ to $\lts$.  
Associated with each function $f \in \ltsg$ is the matrix $M_f \in M(\S)$
whose $(x,y)$ entry is $\fb (d(x,y))$, whence the matrix $M_f$ 
corresponds to the linear operator $h \mapsto h * f$ given in 
Section~\ref{sub:drs}. The following
analogue of Lemma~\ref{lem:spherical} is derived from Section~2.4
of~\cite{Ter98}; a published reference is Section~2.3 of~\cite{BCN}.

\begin{lem} \label{lem:scheme}
There exists a basis of real--valued functions 
$\f_0 , \ldots , \f_D$ of $\ltsg$ 
orthogonal under the inner product 
$\langle f , g\rangle = |\S|^{-1} \sum_x f(x) 
\overline{g(x)}$ with the following properties. \\
(1) $\f_0 (x) \equiv 1$. \\
(2) $\f_j (\zero) = 1 = \sup_x |\f_j (x)|$ for all $j$. \\
%%% This follows from
%%% fact that $E_{ij}$'s are representable as $v_i \cdot v_j$,
%%% or see Terwilliger Lemma 2.11.3 (i) with $u_{ij} = q_j (i) / m_j = 
%%% q_j (i) / q_j (j) = E_{ij} / E_{jj}$.
(3) $\f_i \f_j = \sum_{r=0}^D q^r_{ij} \f_r$ for some nonnegative coefficients
   $q^r_{ij}$ with $\sum_r q^r_{ij} = 1$. \\
(4) $\f_i * \f_j = \gamma_j \delta_{ij} \f_j$, where $\gamma_j : = \f_j * \f_j
   (\zero) = |\S|^{-1} \sum_x \f_j (x)^2$. \\
(5) The functions $\f_j$ are eigenfunctions of any convolution
   operator, that is, $M_f \f_j = c \f_j$ for any $f \in \ltsg$. \\
%%(5) For any $M = M_f \in M(\S)$, $M_f \f_i=c\f_i$ for all $i$. \\
%%the functions $\f_0 , \ldots , \f_D$ are an eigenbasis for $M_f$ for $\ltsg$.
(6) For $f \in \ltsg$, we have $f = \sum_{j=0}^D a_j(f) \f_j$, where
   $a_j : = \gamma_j^{-1} |\S|^{-1} \sum_x f(x) \f_j (x)$. \\
(7) For $f,g \in \ltsg$, we have $a_j(f * g)= \gamma_j a_j(f) a_j(g)$. \\
(8) For $f\in \ltsg$ which is positive and nonincreasing, 
$|\langle f , \f_i\rangle|\le  \langle f , \f_1\rangle $
for each $i\ge 1$.
%%\item The eigenvalues of $M_f$ are 
%%$\gamma_0 a_0(f) , \ldots , \gamma_D a_D(f)$.
\end{lem}

If we place the norm $\sum_{j=0}^D |a_j(f)|$ on $\sspj$, essentially all of the
hypotheses in Theorems~\ref{th:gen ii} and~\ref{th:gen i} (to come later) are 
immediate noting that all norms are equivalent on finite dimensional spaces.
If the analogue of~(\ref{eq:xiv}) holds, then letting
$L (g):= |S|^{-1}\sum_{x\in\S}g(x)\f_1(x)$ and both $\op $ and $\rho$ to
be $L(K_J)$, then one can easily show that {\it all} of the 
hypotheses in Theorems~\ref{th:gen ii} and~\ref{th:gen i} hold.
As far as~(\ref{eq:xiv}), it trivially holds for the complete graph where the 
diameter $D$ is equal to 1 and in any case, the reader is left with only one 
condition to check.

\setcounter{equation}{0}

\section{Two Technical Theorems} \label{sec:proofs}
We now state two general results from which 
Theorems~\ref{th:main} and~\ref{th:potts} will follow.

\begin{th} \label{th:gen ii}
Let $\Gamma$ be any tree (with bounded degree).
For the $d$--dimensional Heisenberg model with $d\ge 1$, 
if $J > 0$ and $$\textstyle {br}  (\Gamma) \cdot \op < 1,$$ 
then there is no robust phase transition for the parameter $J$, where $\op$ is 
given in Definition~\ref{defn:op} ($\op$ implicitly depends on $d$).
More generally, if $J>0$ and
if $(\S , G , H)$ is any statistical ensemble with a norm
$|| \cdot ||$ on $\sspj$ 
satisfying~(\ref{eq:xi}),~(\ref{eq:xii}),~(\ref{eq:present iii}) 
and~(\ref{eq:xiii}) and there exists a number $\op \in (0,1)$ satisfying
(\ref{eq:ixix}) and 
$\textstyle {br}  (\Gamma) \cdot \op < 1$, then there is no robust 
phase transition for the parameter $J$.
\end{th}

\begin{th} \label{th:gen i}
Let $\Gamma$ be any tree (with bounded degree).
For the $d$--dimensional Heisenberg model with $d\ge 1$, 
if $J > 0$ and $$\textstyle {br}  (\Gamma) \cdot \op > 1,$$
then there is a robust phase transition for the parameter $J$, where
$\op$ is as above. More generally, if $J>0$ and 
if $(\S , G , H)$ is any statistical ensemble with a norm
$|| \cdot ||$ on $\sspj$
satisfying~(\ref{eq:xi}),~(\ref{eq:present iii}),
(\ref{eq:present iii'}),~(\ref{eq:xiii}) and~(\ref{eq:xiiii}), 
and if $L$ is a linear
functional on $\sspj$ which vanishes
on the constants and
satisfies~(\ref{eq:(a)}),~(\ref{eq:(c)}) and~(\ref{eq:(b)})
for a constant $\rho > 0$, then 
$\textstyle {br}  (\Gamma) \cdot \rho > 1$ implies
a robust phase transition for the parameter $J$.
\end{th}

To prove these results, 
we begin with a purely geometric lemma on the existence of 
cutsets of uniformly small content below the branching number.

\begin{lem} \label{lem:globalcut}
Assume that $\textstyle {br} (\Gamma) < d$. 
Then for all $\eps >0$, there exists a cutset $C$ such that
$$
\sum_{x\in C}({1 \over d})^{|x|} \le \eps
$$
and
for all $v\in C^i\cup C$,
\be\lab{eqn:goodcut}
\sum_{x\in C\cap \Gamma (v)}({1 \over d})^{|x|-|v|} \le 1.
\ee
\end{lem}

\noindent{\bf Proof.}
Since $\textstyle {br} (\Gamma) < d$, for any given $\eps >0$, there exists 
a cutset $C$ such that
$$
\sum_{x\in C}({1 \over d})^{|x|} \le \eps.
$$
We can assume that $C$ is a minimal cutset with this property with respect
to the partial order $C_1 \preceq C_2$ if for all $v\in C_1$, there exists
$w\in C_2$ such that $v\le w$. We claim that this cutset
satisfies~(\ref{eqn:goodcut}). If this property failed for some $v$,
we let $C'$ be the modified cutset obtained by replacing $C\cap \Gamma (v)$
by $v$ (and leaving $C\cap \Gamma^c_v$ unchanged).
As~(\ref{eqn:goodcut}) clearly holds for $w\in C$, we must
have that $v\not\in C$ in which case $C'\neq C$.
We then have
\begin{eqnarray*}
\sum_{x\in C'}({1 \over d})^{|x|} & = & 
   \sum_{x\in C\cap \Gamma (v)^c}({1 \over d})^{|x|}+({1 \over d})^{|v|} \\[1ex]
& < & \sum_{x\in C\cap \Gamma (v)^c}({1 \over d})^{|x|}+
   ({1 \over d})^{|v|} \sum_{x\in C\cap \Gamma (v)}({1 \over d})^{|v-x|} \\[1ex]
& = & \sum_{x\in C\cap \Gamma (v)^c}({1 \over d})^{|x|}+
   \sum_{x\in C\cap \Gamma (v)}({1 \over d})^{|x|} \\[1ex]
& = & \sum_{x\in C}({1 \over d})^{|x|} \\[1ex]
& \le & \eps,
\end{eqnarray*}
contradicting the minimality of $C$ since clearly $C'\preceq C$.
$\qed$

We now proceed with the proofs of Theorems~\ref{th:gen ii} and~\ref{th:gen i}.

\noindent{\bf Proof of Theorem \ref{th:gen ii}.}  
Since in Section~\ref{sub:sphere}
the Heisenberg models have been shown to satisfy all of
the more general hypotheses of this theorem, we need only prove the last
statement of the theorem where we have a given $J>0$, a given $\|\,\,\|$
on $\sspj$ and a given $\op$ satisfying the required conditions.
By~(\ref{eq:xi}), for any $\eps > 0$, there is
an $\eps_0 > 0$ such that for all $k \leq B$ and all $h_1 , \ldots , 
h_k \in \PP_+(J)$ with $\|h_i - 1\| \leq \eps_0$ for all $i$, we have that
\begin{equation} \label{eq:star2}
\| \poi_k (h_1 , \ldots , h_k) - 1 \| \leq (1 + \eps) \sum_{i=1}^k
   \|h_i - 1\| \, .
\end{equation}
Choose $\eps > 0$ so that 
$(1 + \eps)^{-1} > \textstyle {br}  (\Gamma) \cdot \op$ and 
choose $\eps_0$ as above.
By~(\ref{eq:xii}), we can choose $J'>0$ small enough 
so that $\| K_{J'} -1 \| \leq \eps_0 \op$.
Use Lemma~\ref{lem:globalcut} to choose
a sequence of cutsets $\{ C_n \}$ for which 
$$\lim_{n \rightarrow \infty} \sum_{x \in C_n} [(1 + \eps) 
   \op ]^{|x|} = 0$$
and for all $n$ and all $v \in C_n^i \cup C_n$,
\begin{equation} \label{eq:star43}
\sum_{x \in C_n \cap \Gamma (v)} 
[(1 + \eps) \op ]^{|x| - |v|} \leq 1.
\end{equation}
We now show by induction that for all $n$ and all $v \in C_n^i$, 
\begin{equation} \label{eq:ind}
\|f^{J' ,J, +}_{C_n , v} - 1\| \leq \eps_0 \sum_{x \in C_n \cap \Gamma (v)} 
   [(1 + \eps) \op ]^{|x| - |v|} \, .
\end{equation}
Indeed, from Lemma~\ref{lem:rec}, letting $w_1 , \ldots , w_k$
be the children of $v$,
$$\|f^{J',J , +}_{C_n , v} - 1\| 
= \| \poi (\calk_{J_1''} f^{J',J , +}_{C_n , w_1},
   \ldots , \calk_{J_k''} f^{J',J , +}_{C_n , w_k}) - 1 \| \, $$
where $J_i''$ is $J$ if $w_i\in C_n^i$ and $J'$ otherwise.
When $w_i \in C_n$, the choice of $J'$ guarantees that 
$\|\calk_{J_i''} f^{J',J , +}_{C_n , w_i} - 1\| \leq \eps_0\op\leq \eps_0$, 
while when $w_i \notin C_n$, 
the induction hypothesis together with~(\ref{eq:star43}) guarantees that
$\|f^{J',J , +}_{C_n , w_i}    - 1\| \le \eps_0$ which implies that
$\|\calk_{J_i''} f^{J',J , +}_{C_n , w_i}    - 1\| \le \eps_0$ 
by~(\ref{eq:xiii}). Hence, from~(\ref{eq:star2}), 
$$\|f^{J',J , +}_{C_n , v} - 1\| \leq 
   (1 + \eps) \sum_{w_i \in C_n} \| \KJP f^{J',J , +}_{C_n , w_i} - 1 \| 
   + (1 + \eps) \sum_{w_i \notin C_n} 
\|\K_J f^{J',J , +}_{C_n , w_i}    - 1\|.$$

The summands in the first sum are at most 
$\eps_0 \op$ while those in the second sum are by~(\ref{eq:ixix}) at 
most $\op \|f^{J',J , +}_{C_n , w_i}    - 1\| $.
Therefore using the induction hypothesis on the second term, we obtain
\begin{eqnarray*}
\|f^{J',J , +}_{C_n , v} - 1\| 
& \leq & \sum_{i=1}^k (1 + \eps) \eps_0 \op \sum_{x \in C_n \cap 
   \Gamma (w_i)} \left [ (1 + \eps ) \op \right ]^{|x| - |w_i|} \\[2ex]
& = & \eps_0 \sum_{x \in C_n \cap \Gamma (v)} \left [ (1 + \eps) 
   \op \right ]^{|x| - |v|} ,
\end{eqnarray*}
completing the induction. 
Finally, the theorem follows by taking $v = o$, letting $n \rightarrow
\infty$, and using~(\ref{eq:present iii}). $\qed$

For the proof of Theorem~\ref{th:gen i}, it is easiest to isolate the 
following two lemmas.

\begin{lem} \label{lem:first}
Under the more general hypotheses of Theorem~\ref{th:gen i}
(with a given $J >0$, a given $\|\,\,\|$ on $\sspj$, a given $L$ and 
a given $\rho$ satisfying the required conditions), 
for all $\alpha>0$,
there exists $\beta>0$ so that if $h_1,\ldots, h_k\in \PP_+(J)$ with
$k\le B$ and $\|h_i-1\| < \beta$ for each $i$, then
$$
L \left [ (\poi_k (\KJ h_1 , \ldots , \KJ h_k)) - 1 \right ] 
   \geq {1 \over 1 + \alpha} \sum_{i=1}^k L (\KJ h_i - 1) 
$$
\end{lem}

\noindent{\bf Proof.}
In~(\ref{eq:xi}), choose $\beta<1$ so that 
$$o(h)\le h \left(1-{1 \over (1+\alpha)}\right) {c_4 \over c_3} $$
for all $h\in (0,\beta)$, with $c_3$ and $c_4$ as in~(\ref{eq:(c)})
and~(\ref{eq:(b)}). If $h_1,\ldots, h_k\in \PP_+(J)$ are such that
$\|h_i-1\| < \beta$, then $\|\KJ h_i-1\| < \beta$ by~(\ref{eq:xiii}).
We can now write
\begin{equation} \label{eq:U1}
\poi_k(\KJ h_1,\ldots,\KJ h_k)-1-
{1 \over (1+\alpha)}\sum_{i=1}^k (\KJ h_i-1)
\end{equation}
as 
\begin{equation} \label{eq:U2}
\left(1-{1 \over (1+\alpha)}\right)\sum_{i=1}^k (\KJ h_i-1) +U
\end{equation}
where by assumption,
\begin{eqnarray} \label{eq:new03}
\|U\| & \le & o(\max_i \|\KJ h_i-1 \|) \\[2ex]
& \le & \left(1-{1 \over (1+\alpha)}\right) {c_4 \over c_3} 
   \max_i \|\KJ h_i-1 \| \nonumber \\[2ex]
& \le & \left(1-{1 \over (1+\alpha)}\right){c_4 \over c_3} 
   \sum_{i=1}^k \|\KJ h_i-1 \|. \nonumber
\end{eqnarray}
Letting $a$ be the quantity~(\ref{eq:U1}), we see that
\begin{eqnarray*}
L(a) & = & L \left [ \left ( 1 - {1 \over (1+\alpha)} \right )
   \sum_{i=1}^k (\KJ h_i-1) \right ] + L(U) \\
& \ge & \left ( 1 - {1 \over (1+\alpha)} \right ) c_4 \sum_{i=1}^k 
   \|\KJ h_i-1\| -c_3 \|U\| \\
& \geq & 0
\end{eqnarray*}
by~(\ref{eq:(c)}),~(\ref{eq:(b)}) and~(\ref{eq:new03}), 
which is the conclusion of the lemma.  $\qed$

The next lemma tells us that in ``one step'', we can't move from being ``far
away'' from uniform to being ``very close'' to uniform.

\begin{lem} \label{lem:second}
Under the more general hypotheses of Theorem~\ref{th:gen i}
(with a given $J >0$, a given $\|\,\,\|$ on $\sspj$, a given $L$ and 
a given $\rho$ satisfying the required conditions), 
for all $\beta>0$ and $J'\in (0,J]$, there exists a $\gamma<\beta$ such that if
$\| \poi_k ( \calk_{J_1''} h_1 , \ldots , \calk_{J_k''} h_k) - 1 \| <\gamma$ 
with $h_1,\ldots, h_k\in \PP_+(J)\cup \{\p_{\zero}\}$ and $k\le B$
and with $J_i''$ being $J$ if $h_i\in \PP_+(J)$ and $J'$ if $h_i=\p_{\zero}$,
then each $h_i$ is not $\p_{\zero}$ and $\sum_{i=1}^k \| h_i - 1\| < \beta$.
\end{lem}

\noindent{\bf Proof.}  Choose $\gamma\in (0,\min\{\beta,1/c_1\})$ so that
$$
{2c_1 c_3 B\gamma  \over\rho c_2 c_4 (1-c_1\gamma)} < \beta
$$
and 
$$
\min\{||K_J-1||,||K_{J'}-1||\} > {2c_1 \gamma  \over (1-c_1\gamma)c_2}
$$
where $c_1,c_2,\rho,c_3$ and $c_4$ come 
from~(\ref{eq:present iii}),~(\ref{eq:present iii'}),~(\ref{eq:(a)}),~(\ref{eq:(c)}) and~(\ref{eq:(b)}) respectively.
We first show that if $h_1,\ldots,h_k\in\PP_+(J)$, with $k\le B$, then
$\| \poi_k (h_1 , \ldots , h_k) - 1 \| <\gamma< 1/c_1$ implies that for all 
$i$
$$
\| h_i  - 1 \| <{2c_1\gamma \over (1-c_1\gamma)c_2}.
$$
[Proof:
$$
||h_i-1||\le c_2^{-1}||h_i-1||_\infty
\le c_2^{-1}\left({\max h_i\over \min h_i}-1\right)
$$
$$
\le c_2^{-1}\left({\max \prod_i h_i\over \min \prod_i h_i}-1\right)
=  c_2^{-1}\left({\max \poi_k (h_1 , \ldots , h_k)\over\min \poi_k 
(h_1 , \ldots , h_k) }-1\right)
$$
where the second inequality is straightforward and the
third inequality comes from~(\ref{eq:xiiii}).
Next, $\| \poi_k (h_1 , \ldots , h_k) - 1 \| <\gamma< 1/c_1$ implies 
$|| \poi_k (h_1 , \ldots , h_k) - 1 ||_\infty\le c_1\gamma$  
which implies the last expression is at most 
$$
c_2^{-1}\left({1+c_1\gamma\over 1-c_1\gamma}-1\right)=
c_2^{-1}{2c_1\gamma\over 1-c_1\gamma}.]
$$
It follows that if
$\| \poi_k ( \calk_{J_1''} h_1 , \ldots , \calk_{J_k''} h_k) - 1 \| <\gamma$,
then 
$$
\| \calk_{J_i''} h_i-1 \| <{2c_1 \gamma \over (1-c_1\gamma)c_2}
$$ 
for each $i$ which implies that $h_i\in\PP_+(J)$ 
(as opposed to being $\p_{\zero}$). Hence $J_i''$ is $J$ for all $i$.

Now from~(\ref{eq:(a)})--(\ref{eq:(b)}) we have 
$$||\K_{J} h_i - 1|| \geq {\rho c_4 ||h_i - 1|| \over c_3}$$
and we obtain the conclusion of the lemma.
$\qed$

\noindent{\bf Proof of Theorem~\ref{th:gen i}.}  
Since in Section~\ref{sub:sphere}
the Heisenberg models have been shown to satisfy all of
the more general hypotheses of this theorem, we need only prove the last
statement of the theorem, where we have
a given $J >0$, a given $\|\,\,\|$ on $\sspj$, a given $L$ and 
a given $\rho$ satisfying the required conditions.
Choose an $\alpha > 0$ so that 
$\textstyle {br}  (\Gamma) \cdot \rho > 1 + \alpha$.
Choosing $\beta$ from Lemma \ref{lem:first}, we have, under
our assumptions, that for all $h_1 , \ldots , h_k \in \PP_+(J)$ with
$k\le B$ and $\|h_i-1\| < \beta$ for each $i$, 
\begin{equation} \label{eq:star4old}
L \left [ (\poi_k (\KJ h_1 , \ldots , \KJ h_k)) - 1 \right ] 
   \geq {1 \over 1 + \alpha} \sum_{i=1}^k L (\KJ h_i - 1) 
   \geq {\rho \over 1 + \alpha} \sum_{i=1}^k L (h_i - 1). 
\end{equation}
Now, if there is no robust phase transition, then by~(\ref{eq:present iii'})
there must exist $J'\in (0,J]$ 
and a sequence of cutsets $\{ C_n \}$ going to infinity
such that $\lim_{n \rightarrow \infty} \|f^{J',J , +}_{C_n , o} - 1\| = 0$.
Using Lemma \ref{lem:second}, choose $\gamma<\beta$
corresponding to $\beta$ and $J'$. Next, by our choice of $\alpha$, we have
$$I := \inf_C \sum_{x \in C} \left ( {\rho \over 1 + \alpha} \right)^{|x|} 
   > 0$$
where the infimum is over all cutsets.
We now choose $n$ so that 
$$
\|f^{J',J , +}_{C_n , o} - 1\| < \min \{ \gamma , {c_4 \gamma I \over c_3} \}.
$$
where $c_3$ and $c_4$ come from~(\ref{eq:(c)}) and~(\ref{eq:(b)}) respectively.
We then define $\Gamma'$ to be the component of the set
$$\{ v \in C_n^i : \|f^{J',J , +}_{C_n , v} - 1\| < \gamma \}$$
that contains $o$
and let $C$ be the exterior boundary of $\Gamma'$
(that is, the set of $x \notin \Gamma'$ neighboring some $y \in \Gamma'$).  
By the choice of $\gamma$, $C \subseteq C_n^i$ and 
for each $v\in C^i\cup C$, the density $f^{J',J , +}_{C_n , v}$ is in
$$ \PP_+(J) \cap \{f: \|f-1\| < \beta\}.  $$
Using~(\ref{eq:star4old}) and induction, we see that
$$L (f^{J',J , +}_{C_n , o} - 1) \geq \sum_{x \in C} \left ( 
   {\rho \over 1 + \alpha} \right )^{|x|} L (f^{J',J , +}_{C_n , x} - 1) .$$
By definition of $\Gamma' , C$ and $I$ and the fact that $L (f-1) \geq c_4 
\|f - 1\|$ on $\PP_+(J)$, we see that 
$$L (f^{J',J , +}_{C_n , o} - 1) \geq c_4 \gamma I .$$
Hence
$$\|f^{J',J , +}_{C_n , o} - 1\| \geq {c_4 \over c_3} \gamma I .$$
This contradicts the choice of $n$, proving that there is indeed
a robust phase transition.
$\qed$

\setcounter{equation}{0}

\section{Analysis of specific models} \label{sec:anal}

\subsection{Heisenberg models} \label{sub:spherical}

For the Heisenberg models, recall that $\S = S^d$, $d \geq 1$,
and $H(x,y) = - x \cdot y$.  The operator $\KJ$ is convolution
with the function $K_J (x) = c e^{ J x \cdot \zero}$, where $c$ is 
a normalizing constant.  

\noindent{\bf Proof of Theorem~\protect{\ref{th:main}}.}  
A change of variables shows that $L(K_J)=\calrdJ$ and so the result follows
from Theorems~\ref{th:gen ii} and~\ref{th:gen i}. $\qed$

For the rotor model, we now prove the equivalence of SB and SB+.  

\noindent{\bf Proof of Proposition~\protect{\ref{pr:rotor equiv}}.}   
We have already seen the representation 
$$f = \sum_{n \geq 0} a_n (f) \f_n, $$
for functions $f \in \ltsg$.  In the case of the rotor model, where
$\S = S^1$ and we take $\zero$ to be $(1,0)$, 
the space $\ltsg$ is the space of even functions of
$\theta \in [-\pi , \pi]$ and $\f_n = \cos (n \theta)$.  We now
turn to the full Fourier decomposition $f = \sum_{n \in \Z} b_n (f) 
e^{i n \theta}$, where $b_n (f) = \int_0^{2\pi} f(\theta) e^{-i n \theta} \,
d\theta/ (2 \pi)$.  

Let $C$ be any cutset and $\p$ be a set of boundary conditions on $C$.
Let $\J$ be any set of interaction strengths.  It suffices to show that 
$$||f^{\J , \p}_{C,w} - 1||_\infty \leq ||f^{\J , +}_{C,w}-1||_\infty$$
for all $w\in C^i$.
%% This implies that one gets SB+ for any cutsets but we don't write this.
For $v \in C$ and $n\in \Z$, let $x_{v,n} = b_n (K_{\J(x) , \p (v)})$
where $e$ is the edge from $v$ to its parent.  

\noindent{\em Claim}: For all $y\in C^i$, the
Fourier coefficients $\int_0^{2\pi} e^{i n \theta}
\, d\mu^{\J , \p}_{C,y}(\theta)$, which we denote by 
$\{ u_{y , n} : n \in \Z \}$, are sums of monomials in 
$\{ x_{v,n} \}_{v\in C, n\in\Z}$ with
nonnegative coefficients.  {\it Proof:} Let $w \in C^i$ have children
$w_1 , \ldots , w_r \in C^i$ and $w_{r+1} , \ldots , w_k \in C$.
Then the Fourier coefficients $\{ u_{w,n} : n \in \Z \}$ 
are the convolution of the $k - r$ series $\{ x_{v , n} : n \in \Z \}$ 
as $v$ ranges over $w_{r+1} , \ldots , w_k$, also convolved with the
series $\{ b_n (K_{\J (\overline{wv})}) u_{v,n} : n \in \Z \}$ as $v$
ranges over $w_1 , \ldots , w_r$.  Since $b_n (K_J) \geq 0$, this
establishes the claim via induction and the fundamental recursion.

Now write $x_{v,n}^+$ for the Fourier coefficients $b_n (K_{\J (e)})$ where
$e$ is as before. 
Since $K_{J , e^{i \alpha}} (x) = K_J (e^{-i \alpha} x)$, it follows that
$$|x_{v,n}| = |x_{v,n}^+| .$$
But $x_{v,n}^+$ is real because $K_J$ is even, and has been 
shown to be nonnegative.  Thus
$$|x_{v,n}| = x_{v,n}^+ ,$$
and it follows from the claim that each $u_{w,n}$ has modulus
bounded above by the corresponding $u_{w,n}^+$ when plus boundary
conditions are taken.  Hence
$$
||f^{\J , \p}_{C,w} - 1||_\infty \le ||f^{\J , \p}_{C,w} - 1||_A
\leq \sum_{n \neq 0} |u_{w,n}| 
\leq \sum_{n \neq 0} u_{w,n}^+ = ||f^{\J , +}_{C,w} - 1||_A
= ||f^{\J , +}_{C,w} - 1||_\infty,$$
proving the lemma.   $\qed$

\noindent{\em Remark:} Although we have used special properties
of the Fourier decomposition on $L^2 (S^1)$, there exist similar
decompositions for $S^d$.  We believe that a parallel argument
can probably be constructed, bounding the modulus of the sum of 
the coefficients of spherical harmonics of a given order by the
coefficients one obtains for the analogous monomials in the values 
$a_n (K_{\J (x)})$, whose coefficients are necessarily nonnegative
by the nonnegativity of the connection coefficients $q^r_{ij}$.   
Thus we are led to state:

\begin{pblm} \label{pblm:all spheres}
Prove a version of Proposition~\ref{pr:rotor equiv} for general 
Heisenberg models on trees.
\end{pblm}

\subsection{The Potts model} \label{sub:potts}

\noindent{\bf Proof of Theorem}~\ref{th:potts}.
We will obtain this result from Theorems~\ref{th:gen ii} and~\ref{th:gen i}.
For (i), letting $||\,\,||$ be the $L_\infty$ norm on $\sspj$ and
$\op=\alpha_J$, all of the hypotheses in Theorem~\ref{th:gen ii}
except~(\ref{eq:ixix}) are clear. The function $K_J$ is given by 
$$K_J (x) = c \exp (J (2 \delta_{x,0} - 1))$$
where $c = (e^J + (q-1) e^{-J})^{-1}$.  The operator $\KJ$ is linear and 
$$\KJ \delta_j = c e^J \delta_j + \sum_{i \neq j} c e^{-J} \delta_i \, .$$
Hence in the basis $\delta_0 , \ldots , \delta_{q-1}$, the
matrix representation of $\KJ$  is $c (e^J - e^{-J}) I + c e^{-J} M$
where $M$ is the matrix of all ones.  On the orthogonal complement
of the constant functions, $\KJ$ is $c (e^J - e^{-J}) I$, 
and~(\ref{eq:ixix}) follows, proving (i) by an application of
Theorem~\ref{th:gen ii}.

For (ii), let $||\,\,||$ be the same as above, $\rho=\alpha_J$ and
$L(h)=h(0)-h(1)$. It is then immediate to check that 
all of the hypotheses in Theorem~\ref{th:gen i} hold and we may conclude
(ii) by an application of Theorem~\ref{th:gen i}. $\qed$  

\setcounter{equation}{0}

\section{Proof of Theorem~\protect{\ref{th:0hd}}.} \label{sec:zero}

By Proposition~\ref{prop:SB=} and the fact that any subtree of
a tree with branching number 1 also has branching number 1, 
it suffices to show:
\begin{quote}
For for any $\Gamma$ with $\textstyle {br} (\Gamma) = 1$, 
and any bounded $\J$, there is
a sequence of cutsets $\{ C_n \}$ such that for any sequence 
$\{ \delta_n \}$ of boundary conditions on $\{ C_n \}$, 
$$
\lim_{n\to\infty}\| f^{\J , \delta_n}_{C_n , o} - 1\|_\infty = 0 .
$$
\end{quote}

It is convenient to work with a different measure of size, the {\em Max/Min} 
measure, defined as follows. (This arose already in the proof of
Lemma~\ref{lem:second}.) For any continuous strictly positive 
function $f$ on $\S$, let
$$\|f\|_M := {\max_{x \in \S} f(x) \over \min_{x \in \S} f(x)} \, .$$
It is immediate to see:

\begin{lem} \label{lem:mmequiv}
For any sequence $\{ h_n \}$ of continuous
probability densities, $\|h_n - 1\|_\infty \rightarrow 0$ if and only if 
$\log \|h_n\|_M \rightarrow 0$.
\end{lem}

Next, we examine the effect of $\KJ$ on $\|f\|_M$.  
\begin{lem} \label{lem:unifmm}
For any statistical ensemble $(\S , G , H)$,
any $J_{\rm max}$ and any $T > 0$ there 
is an $\eps > 0$ such that for any continuous strictly 
positive function $f$ with 
$\|f\|_M \leq T$, and any $J \leq J_{\rm max}$, 
$$\log \|\KJ f\|_M \leq (1 - \eps) \log \| f \|_M \, .$$
\end{lem}

\noindent{\bf Proof.}  Fix $H, J$ and $f$ and assume without loss 
of generality that $\int f \, \dx = 1$ since the
{\em Max/Min} measure is unaffected by multiplicative constants.   
Let $[a,b]$ be the smallest closed interval containing the range of $f$
and $[c,d]$ contain the range of $K_J$ with $a,c >0$.
Since $f$ is a probability density,
$a < 1 < b$ (we rule out the trivial case $f \equiv 1$).  Since 
$K_J = c + (1-c) g$ for some probability density $g$, it follows that 
for any $x \in \S$, 
$$c + (1-c) a \leq \KJ f(x) \leq c + (1-c) b .$$
As $J$ varies over $[0 , J_{\rm max}]$, $\min_x K_J (x)$ is bounded 
below by some $c_0 > 0$, so for all such $J$,
$$c_0 + (1-c_0) a \leq \KJ f(x) \leq c_0 + (1-c_0) b$$
and so
$$\| \KJ f \|_M \leq {c_0 + (1 - c_0) b \over c_0 + (1 - c_0) a} \, .$$
Setting $R = \|f\|_M - 1$, we have $b = (1+R) a$ and so 
$$\| \KJ f\|_M \leq {c_0 + (1 - c_0) (1 + R) a \over c_0 + (1 - c_0) a}
   = 1 + R {(1 - c_0) a \over c_0 + (1 - c_0) a}
   \leq 1 + R (1 - c_0) \, .$$
Thus 
\begin{equation} \label{eq:u}
\|\KJ f\|_M \leq 1 + (1 - c_0) \left ( \|f\|_M - 1 \right ) \, .
\end{equation}
The function $\log (1 + (1 - c_0) u) / \log (1 + u)$ is bounded above
by some $1 - \eps < 1$ as $u$ varies over $(0 , T-1]$, and setting
$u = \|f\|_M - 1$ in~(\ref{eq:u}) gives
$$\log \| \KJ f\|_M \leq \log (1 + (1 - c_0) (\|f\|_M - 1)) \leq 
   (1 - \eps) \log \|f\|_M \, ,$$
proving the lemma.  
$\qed$

Proceeding with the proof of Theorem~\ref{th:0hd}, let $C$ be a cutset
with no vertices in the first generation,
$$\partial C = \{ v \in C^i : \exists w \in C~{\rm with }\,\, v \to w \}, $$
and $\p$ be defined on $C$. Clearly, for continuous strictly positive 
functions $h_1 , \ldots , h_k$,
$$
\| \poi (h_1 , \ldots , h_k) \|_M \leq \prod_{i=1}^k \| h_i \|_M \, .
$$
We have also previously seen (Lemma~\ref{lem:unifbd})
that all densities that arise are uniformly
bounded away from 0 and $\infty$ and hence there is a uniform bound on the
$\| \,\, \|_M $ that arise. We can therefore choose $\eps$ from
Lemma~\ref{lem:unifmm}. Next for any $v\in C^i\setminus \partial C$,
applying the fundamental recursion gives
\begin{eqnarray*}
\log \|f^{\J , \delta}_{C,v} \|_M & = & \log \| \poi (\K_{\J (\overline{vw_1})} 
   f^{\J , \delta}_{C , w_1} , \ldots , \K_{\J (\overline{vw_k})} 
   f^{\J , \delta}_{C , w_k} \|_M \\[2ex]
& \leq & \sum_{i=1}^k \log \| \K_{\J (\overline{vw_i})} 
   f^{\J , \delta}_{C , w_i} \|_M \\[2ex]
& \leq & \sum_{i=1}^k (1 - \eps) \log \| f^{\J , \delta}_{C , w_i} \|_M.
\end{eqnarray*}
Working backwards, we find that for any cutset $C$,
$$\log \|f^{\J , \delta}_{C , o}\|_M \leq \sum_{w \in \partial C}
   (1 - \eps)^{|w|} \log \|f^{\J , \delta}_{C , w}\|_M \, .$$
Since $\textstyle {br}  (\Gamma) = 1$ 
one can choose a sequence of cutsets $\{ C_n \}$
such that $\sum_{w \in \partial C_n} (1 - \eps)^{|w|} \rightarrow 0$.
The uniform bound on 
$\|f^{\J , \delta}_{C , w}\|_M$
implies that for any sequence of functions $\p_n$ on $C_n$,
$$
\lim_{n\to\infty}\log \|f^{\J , \delta_n}_{C_n , o} \|_M =0,
$$
which along with Lemma~\ref{lem:mmequiv} proves the theorem.
$\qed$ 

Olle H\"aggstr\"om pointed out to us that this result could also be
obtained using ideas from disagreement percolation.

\setcounter{equation}{0}

\section{Proof of Theorem {\protect{\ref{th:2trees}}}.} \label{sec:potts}

While we assume that $q$ is an integer, the case of nonintegral $q$
can be made sense of via
the random cluster representation, and it is worth noting here that
the break between $q=2$ and $q=3$ happens at $q = 2 + \eps$.  
See~\cite{Ha} for a discussion of the qualitative differences between
the random cluster model on a tree when $q \leq 2$ as opposed to $q > 2$.

\begin{lem} \label{lem:robust}
Assume that all of the hypotheses of Theorem~\ref{th:gen ii} are in force
(in particular,~(\ref{eq:ixix}) and $\textstyle {br}  (\Gamma) \cdot \op < 1$ 
hold and so there is no RPT for the parameter $J$) and in addition that
$\sup_{y\in \S} \|K_{J,y}\| < \infty$ and~(\ref{eq:ixix}) holds for
all $f\in\PP(J)$ (instead of just $\PP_+(J)$).
Then there is a tree $\Gamma'$ with 
$\textstyle {br} (\Gamma') = \textstyle {br} (\Gamma)$ 
such that $\Gamma'$ has no PT for the parameter $J$.
%%, thus no SB or PT for the rotor and Potts models.
\end{lem}

\noindent{\bf Proof.}  We mimic the proof of Theorem~\ref{th:gen ii}. 
Choose $\eps$, $\eps_0$ and
cutsets $\{ C_n \}$ as in the proof of Theorem~\ref{th:gen ii}
where we can assume that the cutsets $\{ C_n \}$ are disjoint.  
Choose an integer $m$ sufficiently large so that the $m$-fold iterated
convolution operator $\KJ^m$ satisfies 
$||\KJ^m \p_y-1|| \leq \eps_0 \op$.
For each increasing sequence $\{ n(k) : k = 1 , 2 , \ldots \}$ 
of integers, define a tree $\Gamma'$ by replacing each edge from an element
of $C_{n(k)}$ to its parent by $m$ edges in series, for all cutsets 
in the sequence $\{ C_{n(k)} \}$.  It is not too great an abuse of 
notation to let $C_n$ denote the cutset of $\Gamma'$ consisting of 
the same vertices as before.  It is now possible to 
establish~(\ref{eq:ind}) for all $v \in D$, where $D$ is the set of 
vertices in $\Gamma'$ that are in $C^i$ and in $\Gamma$ (i.e., are not
in a chain of parallel edges that was added). 
The only adjustment
in the proof is as follows.  Use Lemma~\ref{lem:rec} to represent
$f^{J , +}_{C_n , v}$ in terms of $f^{J,+}_{C_n , w}$ where $w$
are the children of $v$ in $\Gamma$ rather than in $\Gamma'$, i.e., we leap 
the whole chain of $m$ edges at once.  Then the case $w \in C_n$ that 
was handled by the choice of $J'$ is replaced by a case 
$w \in \Gamma' \setminus \Gamma$, which is handled by the choice of $m$.  
In fact, (\ref{eq:ind}) holds when + is replaced by any boundary condition 
as the exact same proof shows.
By choosing $\{ n(k) \}$ sufficiently sparse, we can ensure that 
$\textstyle {br}  (\Gamma') = \textstyle {br}  (\Gamma)$.  
Fixing any such choice of $\{ n(k) \}$,
it follows that there is no phase transition by the above together with
Proposition~\ref{prop:SB=}.
$\qed$

We proceed now with the description of a counterexample.
For $\Gamma_1$, we choose the homogeneous binary tree, where each vertex
has precisely 2 children.  Recall from
Section~\ref{sub:potts} that under + boundary conditions, the functions
$f^{J , +}_{C , v}$ all lie in a one-dimensional set.  The most
convenient parameterization for the segment is by the log-likelihood ratio
of state $\zero$ to the other states.  
Thus the probability measure $a \delta_0 + \sum_{i=1}^{q-1}
((1-a)/(q-1)) \delta_i$ is mapped to the value $\log [(q-1) a / (1-a)]$.
Let $g(v)$ denote the log-likelihood ratio at $v$ under some interaction
strength and boundary conditions.  
The recursion~(\ref{eq:recurse}) of Lemma~\ref{lem:rec} boils down to 
$$g(v) = \sum_{v \rightarrow w} \phi (g(w)) ; \;\;\;\; \phi (z) := \log 
   { p e^z + 1-p \over {1-p \over q-1} e^z + (1 - {1-p \over q-1})} \, ,$$
where 
\begin{equation} \label{eq:p}
p := e^J / (e^J + (q-1) e^{-J}) \, .  
\end{equation}

Taking a Taylor expansion to the second order gives
$$\phi (z) = \left ( p - {1-p \over q-1} \right ) z + {1-p \over 2 (q-1)^2}
   [p(q-1)^2 - (q-1) + (1-p)] z^2 + O(z^3) .$$
To see that the second derivative is positive at 0
for $q > 2$, first take the $q$-derivative of the $z^2$ coefficient which is
$[q+2p-3](1-p)/(2(q-1)^3)$.
The definition of $p$ and the fact that $J > 0$ imply
that $p > 1/q \geq 1/(2(q-1))$. Since $x+1/(x-1) -3 >0$ on $(2,\infty)$
and $2p> 1/(q-1)$,
it follows that the $z^2$ coefficient has a positive
$q$-derivative for $q \ge 2$, and is therefore positive for all $q > 2$.  
(This also implies that for $q\in (2-\p,2)$ for some $\p$,
the function $\phi$ is concave (see~\cite{PP} for a detailed analysis of
the critical case $q=2$).)  

The Taylor expansion gives $\phi'(0) = p - (1-p)/(q-1)$. Note that
$p_0: = (q+1)/(2q)$ satisfies $p_0 - (1 - p_0) / (q-1) = 1/2$.  
The value of $p_0$ is chosen to make $\phi' (0) = 1/2$; by
convexity of $\phi$ near zero, there is an interval 
$I := (p_0 - \eps , p_0)$ such that for $p \in I$, the equation
$\phi (z) = z/2$ has a positive solution, call it $z(p)$. Take
$\eps>0$ so small that $p_0 -\eps > 1/q$.
For any $1>p > 1/q$ there is a unique $J > 0$ such that~(\ref{eq:p})
holds.  If $p \in I$, then $z(p)$ is a fixed point for the function
$2 \phi$ and it is easy to see by induction that under + boundary conditions 
on the binary tree, one will always have $g(v) \geq z(p)$.  Thus
we have shown that $\Gamma_1$ has a phase transition for any
$J$ such that $p \in I$.  

To find $\Gamma_2$, we examine the connection between $p_0$ and 
$\| \KJ \|$ where for the rest of the proof, the operator norm
refers to the $L^\infty$ norm on the orthogonal complement of the constants. 
Observe that
$$p - {1-p \over q-1} = {e^J \over e^J + (q-1) e^{-J}} - {e^{-J} \over
   e^J + (q-1) e^{-J}} = \| \KJ \|$$
by the computation in Section~\ref{sub:potts}.  Thus $p_0$ is
chosen to make $\| \KJ \| = 1/2$ and for any $p \in I$, $\| \KJ \| < 1/2$.
Fix any $J$ so that $p \in I$, and let $\Gamma$ be any tree with 
$$2 = \textstyle {br}  (\Gamma_1) < 
\textstyle {br}  (\Gamma) < \| \KJ \|^{-1}.$$  
Let $\Gamma'$ be as in Lemma~\ref{lem:robust} and set $\Gamma_2 = \Gamma'$.  
Then there is no phase transition on $\Gamma_2$ for the chosen parameters,
and since we have seen there is a phase transition for $\Gamma_1$,
this completes the proof of Theorem~\ref{th:2trees}.  $\qed$

\noindent
{\bf Acknowledgements.}
We thank Richard Askey for discussions and showing us the proof of
Lemma~\ref{lem:incr}, J\"{o}ran Bergh, Yuval Peres
and Paul Terwilliger for discussions, Anton Wakolbinger
for providing us with reference \cite{E} and the referee for a correction
and some suggestions.

\noindent
\begin{tabbing}
enoughs \= fffffffffffffffffffffenoughennnnnnnnnnnnnnn \= \kill
\>  Robin Pemantle               \> Jeffrey E.~Steif \\
\>  Department of Mathematics         \> Department of Mathematics \\
\>  University of Wisconsin-Madison \> Chalmers University of Technology \\
\> Van Vleck Hall               \> S--41296 Gothenburg \\
\>   480 Lincoln Drive                \> Sweden \\
\>   Madison, WI 53706    \> steif@math.chalmers.se \\
\>  pemantle@math.wisc.edu             
\end{tabbing}

\begin{thebibliography}{9}

\bibitem{AVLF}
Adel'son-Vel'skii, G., Veisfeiler, B., Leman, A. and Faradzev, I.
(1969).  Example of a graph without a transitive automorphism group.
{\em Soviet Math. Dokl.} {\bf 10} 440--441.

\bibitem{ACCN} Aizenman, M., Chayes, J. T., Chayes, L. and Newman, C. M.
(1988) Discontinuity of the magnetization in one--dimensional 
$1/|x-y|^2$ Ising and Potts models, {\em J. Stat. Phy.} {\bf 50} 1--40.

\bibitem{Askey74}
Askey, R. (1974).  {\em Orthogonal Polynomials and Special Functions}.
S.I.A.M. Regional conferences in applied mathematics no. 21,
J.W. Arrowsmith, Ltd.: Bristol, England.

\bibitem{BCN}
Brouwer., A., Cohen, A. and Neumaier, A. (1989).  {\em Distance Regular 
Graphs}. Modern Surveys in Mathematics, Ser. 3, Bd. 18.  Springer-Verlag: 
New York.

\bibitem{Big}
Biggs, N. (1993).  {\em Algebraic Graph Theory, 2nd Ed.}  Cambridge 
University Press:Cambridge.  

\bibitem{C} Cassi, D. (1992). 
Phase transition and random walks on graphs: a generalization of the 
Mermin--Wagner theorem to disordered lattices, fractals, and other discrete
structures. {\em Phys. Rev. Lett.} {\bf 68} 3631--3634.

\bibitem{E} Eisele, M. (1994).  
{\em Phase transitions may be absent on graphs with transient random walks}.
{\em Unpublished manuscript.}

%%%\bibitem{Erd}
%%%Erdelyi, A. {\em et al} (1953).  Higher transcendental functions.  
%%%Bateman manuscript project.  McGraw-Hill: New York.
%%%
\bibitem{EKPS} Evans, W., Kenyon, C., Peres, Y. and Schulman, L.J.
(1998). Broadcasting on trees and the Ising model. 
{\em Preprint.}

\bibitem{F}
Furstenberg, H. (1970).
Intersections of Cantor sets and transversality of semigroups.
In {\em Problems in analysis. Sympos. in Honor of Salomon Bochner, 
Princeton Univ. (R. C. Gunning, ed.)} 41--59.
Princeton Univ. Press, Princeton, N.J.

\bibitem{Ge} Georgii, H.-O. (1988).
{\em Gibbs Measures and Phase Transitions}. de Gruyter: New York.

\bibitem{Ha} H\"aggstr\"om, O. (1996).  The random-cluster model on a
homogeneous tree. {\em Probab.\ Theory Related Fields} {\bf 104} 231--253.

%%\bibitem{Ka} Katznelson, Y. (1976).  {\em An Introduction to Harmonic 
%%Analysis}.  Dover: New York.

\bibitem{Lig2} Liggett, T. M. (1996).
Multiple transition points for the contact
process on a binary tree.  {\em Ann. Probab.} {\bf 24} 1675--1710. 

\bibitem{Ly1} Lyons, R. (1989).
The Ising model and percolation on trees and tree-like graphs.
{\em Commun. Math. Phys.} {\bf 125} 337--353.

\bibitem{Ly2}  Lyons, R. (1990).  Random walks and percolation
on trees. {\em Ann. Probab.} {\bf 18} 931--958

\bibitem{MW} Merkl, F., and Wagner, H. (1994).
Recurrent random walks and the absence of continuous symmetry breaking 
on graphs.  {\em J. Stat. Phy.} {\bf 75} 153--165.

\bibitem{MP} Monroe, J. L., and Pearce, P. A. (1979).
Correlation inequalities for vector spin models.
{\em J. Stat. Phy.} {\bf 21} 615--633.

\bibitem{N} Natterer, F. (1986).  {\em The Mathematics of Computerized
Tomography}. John Wiley, Stuttgart.

\bibitem{PS} Patrascioiu A. and Seiler, E. (1992).
Phase structure of two-dimensional spin models and percolation.
{\em J. Stat. Phy.} {\bf 69} 573--595.

\bibitem{Pem} Pemantle, R. (1992). The contact process on trees.
{\em Ann.  Probab.} {\bf 20} 2089--2116.

\bibitem{PP} Pemantle, R., and Peres, Y.,
Recursions on trees and the Ising model, {\em Preprint.}

\bibitem{R} Rainville, E. D. (1960).  {\em Special Functions}.
MacMillan, New York.

\bibitem{Sta}  Stacey, A. (1996). The existence of an intermediate phase for 
the contact process on trees.  {\em Ann. Probab.} {\bf 24} 1711--1726.

\bibitem{Ter98} Terwilliger, P. (1998).  {\em Unpublished lecture notes}.

\end{thebibliography}
\end{document}